\documentclass[11pt]{article}
\usepackage[utf8]{inputenc}
\usepackage{amsmath,amssymb,amsthm,amsfonts}
\usepackage{enumitem}
\usepackage{booktabs}
\usepackage{geometry}
\usepackage{authblk}
\usepackage{appendix}
\usepackage{url}
\usepackage{needspace}
\numberwithin{equation}{section}

\theoremstyle{plain}
\newtheorem{theorem}{Theorem}[section]
\newtheorem{proposition}[theorem]{Proposition}
\newtheorem{lemma}[theorem]{Lemma}
\newtheorem{corollary}[theorem]{Corollary}
\theoremstyle{definition}
\newtheorem{definition}[theorem]{Definition}
\theoremstyle{remark}
\newtheorem{remark}[theorem]{Remark}
\newtheorem{example}[theorem]{Example}

\DeclareMathOperator{\Syz}{Syz}
\DeclareMathOperator{\LM}{LM}
\DeclareMathOperator{\LT}{LT}
\DeclareMathOperator{\LC}{LC}
\DeclareMathOperator{\pd}{pd}
\DeclareMathOperator{\reg}{reg}
\DeclareMathOperator{\Tor}{Tor}

\DeclareMathOperator{\Row}{Row}

\newcommand{\res}[1]{{}_{S}{#1}}         
\newcommand{\cB}{\mathcal{B}}
\newcommand{\cC}{\mathcal{C}}
\newcommand{\N}{\mathbb{N}}

\title{Commutator syzygies and pure resolutions over central polynomial rings with division-ring coefficients}
\author[1]{Hoang Anh Tu\thanks{Email: hoanganhtu.nth.522005@gmail.com.}}
\author[1]{Huynh Viet Khanh\thanks{Email: khanhhv@hcmue.edu.vn.}}
\affil[1]{Department of Mathematics and Informatics, Ho Chi Minh City University of Education, Ho Chi Minh City, Vietnam.}
\date{}

\begin{document}
\maketitle

\begin{abstract}
Let $D$ be a division ring and let $A=D[x,y,z]$ with central indeterminates. For positive integers $p,q$ and $a,b\in D$, put $f=x^p+az^p$ and $g=y^q+bz^q$. We compute their left syzygy module and the minimal graded free resolution of $A/(Af+Ag)$. This quotient has projective dimension two if $ab=ba$ and three otherwise. In the latter case, it has finite length and total Betti numbers $(1,2,2,1)$, with a pure resolution of type $(0,d,3d,4d)$ when $p=q=d$. These tables do not occur over polynomial rings over fields. If $D$ is finite-dimensional over its center $K$, the Betti tables of modules of finite length over $K[x_1,\dots,x_n]$ and left modules of finite length over $D[x_1,\dots,x_n]$ generate the same rational cone for every $n\ge1$. In three variables, the semigroup over $K$ is properly contained in that over $D$ when $D$ is noncommutative.
\end{abstract}



\section{Introduction}

Let $S=k[x_1,\ldots,x_n]$, where $k$ is a field, and suppose that homogeneous polynomials $u,v\in S$ minimally generate an ideal. A syzygy of the pair is a row $(r,s)\in S^2$ such that $ru+sv=0$. If $c$ is a greatest common divisor of $u,v$, their syzygy module is freely generated by $(v/c,-u/c)$. Thus $S/(u,v)$ has projective dimension two and total Betti numbers $(1,2,1)$ (Lemma~\ref{lem:two-generated}). If $u,v$ have relatively prime leading monomials for a monomial order, then $\{u,v\}$ is a Gr\"obner basis of $(u,v)$ \cite[Chapter~2, \S9, Theorem~3, p.~105 and Proposition~4, p.~106]{CLOIVA}. We ask how these conclusions change when the indeterminates remain central and the coefficients lie in a division ring.

A related question concerns which Betti tables can occur over such rings. The diagram $D_1$ in \cite[\S1, equation~(1), p.~343]{ErmanSemigroup} is pure of type $(0,1,3,4)$ and has total Betti numbers $(1,2,2,1)$. Erman shows that $D_1$ cannot be the Betti table of a module over a polynomial ring over a field. He constructs modules with tables $2D_1$ and $3D_1$, and hence realizes $mD_1$ for every integer $m\ge2$ \cite[proof of Theorem~1.6(1), p.~359]{ErmanSemigroup}.

Erman and Sam realize $D_1$ as the Betti table of the residue field of the graded enveloping algebra of the Heisenberg Lie algebra; see \cite[\S6, equation~(2), p.~296]{ErmanSam}. In that algebra, the generators $X,Y$ have degree one and $XY-YX=[X,Y]\ne0$ has degree two. Their example leads us to ask whether $D_1$ can also occur over a polynomial ring with central indeterminates and coefficients in a division ring.

Let $D$ be a division ring and let $A=D[x,y,z]$, with $x,y,z$ central of degree one. Fix positive integers $p,q$ and $a,b\in D$, and put $\kappa=ba-ab$. When $\kappa\ne0$, the polynomials $x^p+a$ and $y^q+b$ generate the unit left ideal in $D[x,y]$. To obtain nonzero graded quotients, we homogenize in $z$ and set
$$
f=x^p+az^p,\qquad g=y^q+bz^q,\qquad I=Af+Ag,\qquad M=A/I.
$$
For any monomial order with $z\prec y\prec x$, the leading monomials of $f,g$ are $x^p,y^q$, which are relatively prime for every choice of $a,b$. Their commutator is $gf-fg=(ba-ab)z^{p+q}=\kappa z^{p+q}$.

When $\kappa=0$, the pair $\{f,g\}$ is a left Gr\"obner basis, $\operatorname{in}_\prec(I)=Ax^p+Ay^q$, and the left syzygy module is generated by $(g,-f)$. When $\kappa\ne0$, the initial ideal is $Ax^p+Ay^q+Az^{p+q}$. These three generators are minimal, and $f,g$ still minimally generate $I$. Thus relatively prime leading monomials do not suffice for $\{f,g\}$ to be a Gr\"obner basis when the coefficients do not commute.

Theorem~\ref{thm:dichotomy} determines the minimal graded free resolution of $M$ for every $D$ and all $p,q\ge1$, with the maps written as right multiplication by matrices. For $\kappa=0$, the module $M$ has projective dimension two and $\dim_D M=\infty$. For $\kappa\ne0$, it has projective dimension three, $\dim_D M=pq(p+q)$, and total Betti numbers $(1,2,2,1)$. If $p=q=d$, the latter resolution is pure of type $(0,d,3d,4d)$. In particular, the case $p=q=1$ realizes $D_1$ over $D[x,y,z]$ for every noncommutative division ring $D$. Corollary~\ref{prop:notcomm} shows that none of the tables with $\kappa\ne0$ occurs over a polynomial ring over a field. We also compute the syzygies of the monic left Gr\"obner basis by Schreyer's construction and then use the formulas for changing generators to pass to the pair $(f,g)$. These constructions for left PBW rings over a division ring are described in \cite[\S4.2]{GomezTorrecillas}.

The commutator calculation also leads to a comparison with the work of Adams and Loustaunau on two Cauchy--Fueter operators. They compute a minimal resolution over the real polynomial ring in eight variables, using columns of matrix commutators in their Gr\"obner basis. The resolution is pure of type $(0,1,3,4)$ with ranks $(4,8,8,4)$ given in \cite[Lemma~2.1, equations~(5) and~(9), and Corollary~3.6]{AdamsLoustaunau}. For $D=\mathbb H$, Proposition~\ref{rem:quaternion-specialization} identifies their resolution, in the ordinary grading, with the restriction of a free resolution over the quaternion polynomial ring in the eight coordinates. After assigning degrees $p,q$ to the coordinates of the two quaternions, the complex specialized at $f,g$ is isomorphic to the resolution of $M$ in Theorem~\ref{thm:dichotomy}(2) when $\kappa\ne0$.

The distinction between the semigroup of Betti diagrams and its rational cone, considered in \cite[Definitions~1.1--1.2, p.~343]{ErmanSemigroup}, leads us to compare both objects when the coefficient field is replaced by a division ring. Suppose that $K$ is a subfield of the center of $D$ with $\delta=\dim_K D<\infty$. Theorem~\ref{thm:restriction} shows that restriction from $D[x_1,\dots,x_n]$ to $K[x_1,\dots,x_n]$ multiplies the Betti table by $\delta$, and extension of scalars preserves the table. For modules of finite length, Corollary~\ref{cor:betti-cones} therefore identifies the two rational cones for every $n\ge1$. In three variables, the semigroup over $K$ is properly contained in that over $D$ when $D$ is noncommutative. The example realizing $D_1$ proves that this inclusion is strict.

  The paper is organized as follows. Section~\ref{sec:graded} contains the preliminary results on graded left modules and syzygies. In Section~\ref{sec:dichotomy}, we compute the minimal resolutions for commuting and noncommuting coefficients and compare their Betti tables with those over fields. Section~\ref{sec:pure} treats the quaternion specialization, and Section~\ref{sec:center} compares Betti tables under extension and restriction of scalars, with applications to semigroups and rational cones.

For the general results, we write $A=D[x_1,\dots,x_n]=\bigoplus_{\nu\ge0}A_\nu$, where $D$ is a division ring, $n\ge1$, and the indeterminates are central. For $\alpha\in\mathbb N^n$, put $x^\alpha=x_1^{\alpha_1}\cdots x_n^{\alpha_n}$. We use the standard grading, with $A_0=D$ and $\deg x^\alpha=|\alpha|=\alpha_1+\cdots+\alpha_n$. The ring homomorphism $A\to D$ defined by $\sum_\alpha a_\alpha x^\alpha\mapsto a_0$ has kernel $A_+=\bigoplus_{\nu>0}A_\nu$, so $A/A_+\cong D$. For a graded left $A$-module $M$ and $d\in\mathbb Z$, the shift $M(-d)$ is the module $M$ with grading $M(-d)_j=M_{j-d}$. In particular, $1\in A(-d)_d$ generates $A(-d)$ as a left $A$-module. We write elements of finite free left modules as rows and lists of generators as columns.

Appendix~\ref{app:code} describes the computational checks. The Python scripts and accompanying documentation are available on Zenodo~\cite{TuKhanhPython}.

\section{Graded left modules and syzygies}\label{sec:graded}\label{pg:s2}

\subsection{Notation and elementary facts}\label{subsec:notation}
\begin{lemma}\label{lem:monomial-basis}
Every $f\in A$ has a unique expression $f=\sum_\alpha a_\alpha x^\alpha$, where $a_\alpha\in D$ and only finitely many coefficients are nonzero. The monomials form both a left and a right $D$-basis of $A$. Moreover, $x^\beta\in x^\alpha A$ if and only if $x^\beta\in Ax^\alpha$, and these conditions hold if and only if $\alpha_i\le\beta_i$ for all $i$. In that case $x^\beta=x^{\beta-\alpha}x^\alpha$. This defines \emph{divisibility of monomials} on either side.
\end{lemma}
\begin{proof}
For $f\in A$, write $f=\sum_\gamma a_\gamma x^\gamma$ with coefficients written on the left. Since $a_\gamma x^\gamma=x^\gamma a_\gamma$, we also have $f=\sum_\gamma x^\gamma a_\gamma$. If $f=\sum_\gamma x^\gamma b_\gamma$, then
$$
0=\sum_\gamma x^\gamma(a_\gamma-b_\gamma)=\sum_\gamma(a_\gamma-b_\gamma)x^\gamma.
$$
It follows that $a_\gamma=b_\gamma$ for every $\gamma\in\mathbb N^n$. Since $x^\alpha$ is central, $x^\alpha A=Ax^\alpha$. Multiplying both sides of $A=\bigoplus_{\eta\in\mathbb N^n}D x^\eta$ on the right by $x^\alpha$, we obtain $Ax^\alpha=\bigoplus_{\eta\in\mathbb N^n}D x^{\eta+\alpha}$. Thus $x^\beta\in Ax^\alpha$ if and only if $\beta-\alpha\in\mathbb N^n$, or equivalently $\alpha_i\le\beta_i$ for $1\le i\le n$. In this case $x^\beta=x^{\beta-\alpha}x^\alpha$.
\end{proof}

To define the leading monomial of a nonzero polynomial, we use an order on the monomials of $A$ that is preserved by multiplication by a monomial.

\begin{definition}\label{def:momomialorder}
A monomial order $\prec$ is a total order on the set of monomials $\{x^\alpha\mid\alpha\in\mathbb N^n\}$ with the following properties.
\begin{enumerate}[label=\textup{(\arabic*)}]
\item $x^\alpha\prec x^\beta$ implies $x^{\alpha+\gamma}\prec x^{\beta+\gamma}$ for every $\gamma\in\mathbb N^n$;
\item every nonempty set of monomials has a least element.
\end{enumerate}
\end{definition}

\begin{example}
On the monomials of $D[x,y]$, define the lexicographic order by
$$
x^a y^b\prec x^c y^d
\quad\Longleftrightarrow\quad
a<c\quad\text{or}\quad(a=c\text{ and }b<d).
$$
For example, $y^3\prec x\prec xy$. Multiplying both monomials by $x^r y^s$ adds $(r,s)$ to each exponent pair and preserves the comparison. For a nonempty set $T$ of monomials, let $a_0$ be the least exponent of $x$ among its elements. With $b_0=\min\{b\mid x^{a_0}y^b\in T\}$, the least element of $T$ is $x^{a_0}y^{b_0}$.
Comparing total degrees first and using $\prec$ for equal degrees defines another monomial order, the graded lexicographic order $\prec_{\mathrm{grlex}}$, in which $x\prec_{\mathrm{grlex}}y^3$.
\end{example}

Fix a monomial order $\prec$ on the monomials of $A$. For $0\ne f=\sum_\alpha a_\alpha x^\alpha$, let $\LM(f)$ be the largest monomial with nonzero coefficient, let $\LC(f)$ be that coefficient, and put $\LT(f)=\LC(f)\LM(f)$. We call $f$ \emph{monic} if $\LC(f)=1_D$.

\begin{proposition}\label{prop:leading-product}
For nonzero $u,v\in A$ one has
$$
\LM(uv)=\LM(u)\LM(v),\qquad \LC(uv)=\LC(u)\LC(v),\quad\text{and}\quad \LT(uv)=\LT(u)\LT(v).
$$
\end{proposition}
\begin{proof}
Write $u=\sum_\eta a_\eta x^\eta$ and $v=\sum_\theta b_\theta x^\theta$, with $\LM(u)=x^\alpha$ and $\LM(v)=x^\beta$. Since the monomials are central, $uv=\sum_{\eta,\theta}a_\eta b_\theta x^{\eta+\theta}$. For $a_\eta,b_\theta\ne0$, multiplicativity of the monomial order implies $x^{\eta+\theta}\preceq x^{\alpha+\theta}\preceq x^{\alpha+\beta}$, and $x^{\eta+\theta}\prec x^{\alpha+\beta}$ if $\eta\ne\alpha$ or $\theta\ne\beta$. The coefficient of $x^{\alpha+\beta}$ is therefore $a_\alpha b_\beta\ne0$. Hence $\LM(uv)=x^{\alpha+\beta}$, $\LC(uv)=a_\alpha b_\beta$, and $\LT(uv)=a_\alpha b_\beta x^{\alpha+\beta}=\LT(u)\LT(v)$.
\end{proof}

\begin{corollary}\label{cor:elementary}
\begin{enumerate}[label=\textup{(\arabic*)}]
\item The ring $A$ is a domain.
\item One has $\deg(uv)=\deg u+\deg v$ for nonzero $u,v\in A$.
\item One has $A^\times=D^\times$.
\end{enumerate}
\end{corollary}
\begin{proof}
For nonzero $u,v\in A$, Proposition~\ref{prop:leading-product} implies $uv\ne0$, so $A$ is a domain. Write $u=\sum_{i=0}^r u_i$ and $v=\sum_{j=0}^s v_j$, with $u_i\in A_i$, $v_j\in A_j$, and $u_r,v_s\ne0$. Each product $u_iv_j$ lies in $A_{i+j}$. Since $i\leq r$ and $j\leq s$, we have $i+j\leq r+s$, and $i+j=r+s$ if and only if $(i,j)=(r,s)$. Thus $(uv)_{r+s}=u_rv_s$ and $(uv)_k=0$ for $k>r+s$. As $u_rv_s\ne0$, we have $\deg(uv)=r+s=\deg u+\deg v$.

If $u\in A^\times$, we apply \textup{(2)} to obtain $0=\deg(uu^{-1})=\deg u+\deg(u^{-1})$. Since $0\leq\deg u\leq\deg u+\deg(u^{-1})=0$, we have $\deg u=0$ and hence $u\in A_0\setminus\{0\}=D^\times$. Conversely, the inverse of an element of $D^\times$ lies in $D\subseteq A$.
\end{proof}

\begin{remark}\label{lem:unit-leading-regular}\label{lem:noetherian}
\begin{enumerate}[label=\textup{(\arabic*)}]
\item Let $B$ be a ring with identity and $t$ a central indeterminate. If $q=ut^r+\cdots$ with $u\in B^\times$, then the right multiplication map $L_q\colon B[t]\to B[t]$ defined by $L_q(v)=vq$ is injective. 
\item Since $D$ is left noetherian, so is $A$ by \cite[Chapter~1, Theorem~2.9(iv)]{McConnellRobson}.
\end{enumerate}
\end{remark}

\subsection{Minimal graded resolutions}\label{subsec:minimal-resolutions}
In \S\S~\ref{subsec:minimal-resolutions}--\ref{subsec:betti-hilbert}, we assume that modules are finitely generated graded left $A$-modules. For such modules $M$ and $N$, we consider an $A$-linear map $\varphi\colon M\to N$ of degree $0$, so $\varphi(M_j)\subseteq N_j$ for every $j$.

For a module $M$, we first establish the graded isomorphism $D\otimes_A M\cong M/A_+M$. For a graded right $A$-module $L$, the component $(L\otimes_A M)_k$ is generated as an abelian group by the tensors $w\otimes m$ with $w\in L_i$, $m\in M_j$, and $i+j=k$. For $w\in L_i$, $m\in M_j$, and $a\in A_\ell$, the element $(wa)\otimes m-w\otimes(am)$ lies in degree $i+j+\ell$ of $L\otimes_{\mathbb Z}M$. By additivity, these elements generate the kernel of the quotient map $L\otimes_{\mathbb Z}M\to L\otimes_A M$, so $L\otimes_A M$ inherits the grading. Regard $D=A/A_+$ as a graded $(D,A)$-bimodule, and write $\overline a=a+A_+$ for $a\in A$. Consider the maps
$$
\begin{aligned}
&\Phi\colon D\otimes_A M\longrightarrow M/A_+M
\quad\text{defined by}\quad \Phi(\overline a\otimes m)=am+A_+M,\\
&\Psi\colon M/A_+M\longrightarrow D\otimes_A M
\quad\text{defined by}\quad \Psi(m+A_+M)=\overline1\otimes m.
\end{aligned}
$$
The action of $D$ on $M/A_+M$ induces the graded $D$-linear map $\Phi$. For $a\in A_+$, we have $\overline1\otimes am=\overline a\otimes m=0$, so $\Psi$ is well defined. The formulas imply $\Phi\Psi(m+A_+M)=m+A_+M$ and $\Psi\Phi(\overline a\otimes m)=\overline1\otimes am=\overline a\otimes m$. Thus $D\otimes_A M\cong M/A_+M$ as graded left $D$-vector spaces. In particular, $D\otimes_AA(-d)\cong D(-d)$.

A module is \emph{graded free} if it has a homogeneous basis; see \cite[\S1.2.4]{Hazrat}.

For a ring $R$ with identity and $P\in M_{r\times s}(R)$, the map $\rho_P\colon R^{1\times r}\to R^{1\times s}$ defined by $\rho_P(u)=uP$ is left $R$-linear. For a compatible matrix $Q$, we have $\rho_Q\circ\rho_P=\rho_{PQ}$.

\begin{lemma}\label{lem:matrix-degree}
Let $F=\bigoplus_{i=1}^rA(-d_i)$ and $G=\bigoplus_{j=1}^sA(-d'_j)$ be graded free left $A$-modules, and let $P=(p_{ij})\in M_{r\times s}(A)$. The map $\rho_P\colon F\to G$ has degree $0$ if and only if each $p_{ij}$ is zero or homogeneous of degree $d_i-d'_j$.
\end{lemma}
\begin{proof}
Let $e_i$ and $e'_j$ be the canonical basis vectors of $F$ and $G$, of degrees $d_i$ and $d'_j$, respectively. The $j$th coordinate of $\rho_P(e_i)=\sum_j p_{ij}e'_j$ is $p_{ij}$, and $G_{d_i}=\bigoplus_j A_{d_i-d'_j}e'_j$. Thus $\rho_P(e_i)\in G_{d_i}$ precisely when $p_{ij}\in A_{d_i-d'_j}$ for every $j$. If $\rho_P$ has degree $0$, then $e_i\in F_{d_i}$ implies $\rho_P(e_i)\in G_{d_i}$, so $p_{ij}\in A_{d_i-d'_j}$ for all $i,j$. Conversely, suppose $p_{ij}\in A_{d_i-d'_j}$ for all $i,j$. An element $w=\sum_i a_ie_i\in F_k$ has $a_i\in A_{k-d_i}$, so $a_i\rho_P(e_i)\in G_k$. Multiplying both sides of $w=\sum_i a_i e_i$ on the right by $P$, we obtain $\rho_P(w)=\sum_i a_i(e_iP)=\sum_i a_i\rho_P(e_i)\in G_k$. Thus $\rho_P$ has degree $0$.
\end{proof}

Let $F$ and $G$ be graded free left $A$-modules with homogeneous bases $u_1,\dots,u_r$ and $v_1,\dots,v_s$, respectively. Put $d_i=\deg u_i$ and $d'_j=\deg v_j$. The map $\bigoplus_{i=1}^r A(-d_i)\to F$ defined by $(a_1,\dots,a_r)\mapsto\sum_{i=1}^r a_i u_i$ is an isomorphism of degree $0$. The basis of $G$ identifies $G$ with $\bigoplus_{j=1}^s A(-d'_j)$ in the same way. For an $A$-linear map $\varphi\colon F\to G$, let $P=(p_{ij})$ be the matrix defined by $\varphi(u_i)=\sum_{j=1}^s p_{ij}v_j$. Then
$$
\varphi\Bigl(\sum_{i=1}^r a_i u_i\Bigr)
=\sum_{j=1}^s\Bigl(\sum_{i=1}^r a_i p_{ij}\Bigr)v_j,
$$
so the row of coordinates $(a_1,\dots,a_r)$ is sent to $(a_1,\dots,a_r)P$. Thus Lemma~\ref{lem:matrix-degree} applies in any homogeneous bases.

\begin{definition}\label{def:graded-free}
A \emph{graded free resolution} of $M$ is an exact sequence $\cdots\to F_2\xrightarrow{\partial_2}F_1\xrightarrow{\partial_1}F_0\xrightarrow{\varepsilon}M\to0$ of degree-zero maps in which every $F_i$ has a finite homogeneous basis. We put $\partial_0=\varepsilon$ and call the resolution \emph{minimal} if $\partial_i(F_i)\subseteq A_+F_{i-1}$ for all $i\ge1$.
\end{definition}

A generating set of a module $M$ is \emph{minimal} if no proper subset generates $M$.

\begin{proposition}\label{prop:generators}
Let $M$ be a finitely generated graded left $A$-module and put $\overline M=M/A_+M$. The following hold.
\begin{enumerate}[label=\textup{(\arabic*)}]
\item The module $M$ has a finite homogeneous generating set. There is an integer $d_0$ such that $M_j=0$ for all $j<d_0$.
\item \textup{(Graded Nakayama lemma)} If $N$ is a graded submodule of $M$ and $M=N+A_+M$, then $M=N$.
\item Let $m_1,\dots,m_r$ be distinct homogeneous elements of $M$, and put $\overline m_i=m_i+A_+M$. Then
$$
M=\sum_{i=1}^r Am_i
\quad\Longleftrightarrow\quad
\overline M=\sum_{i=1}^r D\overline m_i.
$$
The set $\{m_1,\dots,m_r\}$ is a minimal generating set of $M$ if and only if $(\overline m_i)_{i=1}^r$ is a $D$-basis of $\overline M$.
\item If $S$ is a minimal homogeneous generating set of $M$, then $|S|=\dim_D\overline M$ and $|S\cap M_d|=\dim_D\overline M_d$ for each $d\in\mathbb Z$. Every generating set of $M$ has at least $\dim_D\overline M$ elements. If $F_\bullet\to M$ is a minimal graded free resolution and $i\ge1$, the images under $\partial_i$ of any homogeneous basis of $F_i$ form a minimal homogeneous generating set of $\operatorname{im}\partial_i$.
\end{enumerate}
\end{proposition}
\begin{proof}
\textup{(1)} Choose generators $u_1,\dots,u_s$ of $M$ and write $u_a=\sum_k u_{a,k}$ with $u_{a,k}\in M_k$. The nonzero components $u_{a,k}$ generate $M$, since their $A$-span contains every $u_a$. Thus $M_j=\sum_{a,k}A_{j-k}u_{a,k}$. Only finitely many $u_{a,k}$ are nonzero, so there is an integer $d$ with $k\ge d$ whenever $u_{a,k}\ne0$. For $j<d$, we have $A_{j-k}=0$ for each such $k$, hence $M_j=0$.

\textup{(2)} The graded quotient $Q=M/N$ satisfies $Q=A_+Q$. If $Q\ne0$, let $d$ be its least nonzero degree, which exists by \textup{(1)}. Since $Q_{d-\nu}=0$ for every $\nu>0$, we have
$$
0\ne Q_d=(A_+Q)_d=\sum_{\nu>0}A_\nu Q_{d-\nu}=0.
$$
Thus $Q=0$ and $M=N$. The corresponding forms of the graded Nakayama lemma for right and left modules are stated in \cite[\S3.8]{Hazrat} and \cite[Corollary~2.9.2]{NvO}, respectively.

\textup{(3)} Put $\overline M=M/A_+M$ and write $\overline m_j$ for the class of $m_j$. Since $A$ acts on $\overline M$ through $D$, the image of $N=\sum_j Am_j$ in $\overline M$ is $\sum_j D\overline m_j$. Thus the classes span $\overline M$ if and only if $M=N+A_+M$, which is equivalent to $M=N$ by \textup{(2)}.

Suppose that the classes $\overline m_1,\dots,\overline m_r$ span $\overline M$. Applying the spanning criterion to the elements other than $m_k$, we have
$$
M=\sum_{j\ne k}Am_j
\quad\Longleftrightarrow\quad
\overline M=\sum_{j\ne k}D\overline m_j
\quad\Longleftrightarrow\quad
\overline m_k\in\sum_{j\ne k}D\overline m_j.
$$
If $\sum_j c_j\overline m_j=0$ with $c_k\ne0$, then $\overline m_k=-\sum_{j\ne k}(c_k^{-1}c_j)\overline m_j$. Conversely, expressing $\overline m_k$ as a combination of the other classes yields a relation with coefficient $1$ on $\overline m_k$. Thus no $m_k$ can be omitted precisely when the classes are linearly independent over $D$, hence a basis of $\overline M$.

By \textup{(1)}, the images of a finite homogeneous generating set span $\overline M$, so $\dim_D\overline M<\infty$. Choose a homogeneous $D$-basis of $\overline M$ and lift each basis element to $M$ in the same degree. By the minimality criterion, these lifts form a minimal homogeneous generating set of $M$.

\textup{(4)} By \textup{(3)}, the classes of a minimal homogeneous generating set form a basis of $\overline M$, with $\dim_D\overline M_d$ elements in degree $d$ and $\dim_D\overline M$ in total. The image of any generating set spans $\overline M$ over $D$, so every generating set has at least $\dim_D\overline M$ elements.

Let $F_\bullet$ be a minimal graded free resolution of $M$. For $i\ge1$, put $N_i=\operatorname{im}\partial_i$ and regard $\partial_i$ as a surjection $F_i\to N_i$. By exactness and minimality, $\ker\partial_i=\operatorname{im}\partial_{i+1}\subseteq A_+F_i$. Since $\partial_i(A_+F_i)=A_+N_i$, $\partial_i^{-1}(A_+N_i)=A_+F_i+\ker\partial_i=A_+F_i$. Therefore $\partial_i$ induces the graded isomorphism
$$
F_i/A_+F_i\xrightarrow{\ \sim\ }N_i/A_+N_i
\quad\text{defined by}\quad u+A_+F_i\longmapsto\partial_i(u)+A_+N_i.
$$
If $(e_a)$ is a homogeneous basis of $F_i$, its classes form a $D$-basis of $F_i/A_+F_i$. By \textup{(3)}, the elements $\partial_i(e_a)$ form a minimal homogeneous generating set of $N_i$.
\end{proof}

\begin{proposition}\label{prop:minimal-resolutions}
Let $M$ be a finitely generated graded left $A$-module. The following hold.
\begin{enumerate}[label=\textup{(\arabic*)}]
\item For a graded free resolution $F_\bullet\to M$, let $P_i$ be the matrix of $\partial_i$ in homogeneous bases. The resolution is minimal if and only if $P_i\equiv0\pmod{A_+}$ for every $i\ge1$. For each $i\ge1$, this congruence is equivalent to $D\otimes_A\partial_i=0$, and also to every entry of $P_i$ being a nonunit.
\item Any two minimal graded free resolutions $F_\bullet\to M$ and $G_\bullet\to M$, with augmentations $\varepsilon_F$ and $\varepsilon_G$, are isomorphic as graded resolutions of $M$. If $\alpha_\bullet\colon F_\bullet\to G_\bullet$ is a chain map of internal degree $0$ with $\varepsilon_G\alpha_0=\varepsilon_F$, then $\alpha_\bullet$ is an isomorphism.
\item The module $M$ has a minimal graded free resolution $F_\bullet\to M$ with a finite homogeneous basis in each term. For any such resolution, let $\partial_0$ denote the augmentation. Then
$$
\operatorname{rank}_A F_i=\dim_D(N_i/A_+N_i)\quad\text{for all }i\ge0,
$$
where $N_i=\operatorname{im}\partial_i$. These ranks are independent of the choice of minimal resolution.
\end{enumerate}
\end{proposition}
\begin{proof}
\textup{(1)} Let $F_\bullet\to M$ be a graded free resolution and fix $i\ge1$. Choose homogeneous bases $(e_a)$ of $F_i$ and $(e'_b)$ of $F_{i-1}$, and write $\partial_i(e_a)=\sum_b p_{ab}e'_b$, so $P_i=(p_{ab})$. The inclusion $\partial_i(F_i)\subseteq A_+F_{i-1}$ holds precisely when $p_{ab}\in A_+$ for all $a,b$, that is, when $P_i\equiv0\pmod{A_+}$. In $D\otimes_A F_{i-1}$ we have
$$
(D\otimes_A\partial_i)(\overline1\otimes e_a)
=\sum_b\overline1\otimes p_{ab}e'_b
=\sum_b\overline{p_{ab}}(\overline1\otimes e'_b).
$$
The tensors $\overline1\otimes e_a$ form a $D$-basis of $D\otimes_A F_i$, and the tensors $\overline1\otimes e'_b$ form a $D$-basis of $D\otimes_A F_{i-1}$. Hence $D\otimes_A\partial_i=0$ if and only if $\overline{p_{ab}}=0$ for all $a,b$. By Lemma~\ref{lem:matrix-degree}, each nonzero $p_{ab}$ is homogeneous. Since $A^\times=D^\times$ by Corollary~\ref{cor:elementary}, we have $p_{ab}\notin A^\times$ if and only if $p_{ab}\in A_+$.

\textup{(2)} Since $A$ is left noetherian, finitely generated graded left $A$-modules with maps of degree $0$ form an abelian category. A map $A(-d)\to L$ of degree $0$ is determined by the image of $1$ in $L_d$. For any graded surjection $L\to L'$, the map $L_d\to L'_d$ is surjective, so $A(-d)$ is projective in this category. By Proposition~\ref{prop:generators}\textup{(1)}, every object is a quotient of a finite direct sum of shifts of $A$.

Let $F_\bullet\to M$ be a minimal graded free resolution and let $G_\bullet\to M$ be a graded free resolution. By the comparison theorem \cite[Theorem~2.2.6]{Weibel}, there are chain maps $\alpha\colon F_\bullet\to G_\bullet$ and $\beta\colon G_\bullet\to F_\bullet$ of degree $0$ lifting $\mathrm{id}_M$. For any such pair, $\beta\alpha$ is homotopic to the identity on the unaugmented complex $F_\bullet$. Put $\varphi_i=\beta_i\alpha_i$. The homotopy maps $h_i\colon F_i\to F_{i+1}$ have degree $0$. For each $i\ge0$ they satisfy
$$
\varphi_i-\mathrm{id}_{F_i}=\partial_{i+1}h_i+h_{i-1}\partial_i,
$$
where the second term is omitted for $i=0$.

By minimality, $\partial_{i+1}(F_{i+1})\subseteq A_+F_i$. For $i\ge1$, minimality and $A$-linearity of $h_{i-1}$ imply $h_{i-1}(\partial_i(F_i))\subseteq h_{i-1}(A_+F_{i-1})\subseteq A_+F_i$. Thus $(\varphi_i-\mathrm{id}_{F_i})(F_i)\subseteq A_+F_i$, so $F_i=\varphi_i(F_i)+A_+F_i$. By the graded Nakayama lemma in Proposition~\ref{prop:generators}\textup{(2)}, $\varphi_i$ is surjective. Since $F_i$ is noetherian, the chain $\ker\varphi_i\subseteq\ker\varphi_i^2\subseteq\cdots$ stabilizes. Choose $r\ge1$ with $\ker\varphi_i^r=\ker\varphi_i^{r+1}$. If $u\in\ker\varphi_i$, choose $v\in F_i$ with $\varphi_i^r(v)=u$, using surjectivity. Then $\varphi_i^{r+1}(v)=0$, so $v\in\ker\varphi_i^r$ and $u=0$. Hence each $\varphi_i$ is an isomorphism, and so is $\beta\alpha$. If $G_\bullet$ is also minimal, interchanging $F_\bullet$ and $G_\bullet$ shows that $\alpha\beta$ is an isomorphism. Thus $\alpha$ has left inverse $(\beta\alpha)^{-1}\beta$ and right inverse $\beta(\alpha\beta)^{-1}$, so $\alpha$ is an isomorphism.

\textup{(3)} Set $K_{-1}=M$. Given a finitely generated graded module $K_{i-1}$ for $i\ge0$, choose homogeneous elements $m_{i,1},\dots,m_{i,s_i}$ whose classes form a $D$-basis of $K_{i-1}/A_+K_{i-1}$. Put $F_i=\bigoplus_{a=1}^{s_i}A(-\deg m_{i,a})$, with standard generators $e_a$. Let
$$
\pi_i\colon F_i\longrightarrow K_{i-1}
\quad\text{be defined by}\quad
\pi_i(e_a)=m_{i,a}.
$$
Since $\deg e_a=\deg m_{i,a}$, the map $\pi_i$ has degree $0$. By Proposition~\ref{prop:generators}\textup{(3)}, it is surjective. Its reduction $F_i/A_+F_i\to K_{i-1}/A_+K_{i-1}$ sends the classes of the $e_a$ to the chosen basis, so it is an isomorphism. In particular, $K_i=\ker\pi_i\subseteq A_+F_i$. The submodule $K_i$ is graded and is finitely generated since $F_i$ is noetherian, so the construction continues by induction.

Set $\partial_0=\pi_0$. For $i\ge1$, define $\partial_i\colon F_i\to F_{i-1}$ by composing $\pi_i$ with the inclusion $K_{i-1}\hookrightarrow F_{i-1}$. Since each $\pi_i$ is surjective and the inclusions are injective,
$$
\operatorname{im}\partial_0=M
\quad\text{and}\quad
\operatorname{im}\partial_{i+1}=K_i=\ker\partial_i\subseteq A_+F_i\quad\text{for all }i\ge0.
$$
Thus $F_\bullet\to M$ is a minimal graded free resolution.

For any minimal graded free resolution $F_\bullet\to M$, put $N_i=\operatorname{im}\partial_i$ for $i\ge0$, so $N_0=M$. We have $\partial_i(A_+F_i)=A_+N_i$ and $\ker\partial_i=\operatorname{im}\partial_{i+1}\subseteq A_+F_i$. It follows that the map $F_i/A_+F_i\to N_i/A_+N_i$ induced by $\partial_i$ is an isomorphism. Taking $D$-dimensions, we obtain
$$
\operatorname{rank}_A F_i=\dim_D(F_i/A_+F_i)=\dim_D(N_i/A_+N_i).
$$
By \textup{(2)}, the rank of $F_i$ is independent of the choice of minimal resolution.
\end{proof}

\subsection{Betti numbers and Hilbert series}\label{subsec:betti-hilbert}
For a module $M$, write $H_M(t)=\sum_{j\in\mathbb Z}(\dim_D M_j)t^j$ for its Hilbert series. Each $M_j$ is finite-dimensional over $D$ because $M$ is a graded quotient of a finite direct sum of shifts of $A$ and $\dim_D A_j<\infty$ for every $j$. Let $F_\bullet\to M$ be a minimal graded free resolution with $F_i\cong\bigoplus_j A(-j)^{b_{i,j}}$. Here $i$ denotes homological degree and $j$ internal degree. The grading on $D\otimes_A F_\bullet$ induces the internal grading on
$$
\Tor_i^A(D,M)=H_i(D\otimes_A F_\bullet)\quad\text{for }i\ge0.
$$
By Proposition~\ref{prop:minimal-resolutions}\textup{(2)}, any two minimal resolutions are related by a chain isomorphism of degree $0$ lifting $\mathrm{id}_M$. After tensoring with $D$, this map induces a graded isomorphism on homology.

For the minimal resolution, we obtain that $D\otimes_A\partial_i=0$ for $i\ge1$ by Proposition~\ref{prop:minimal-resolutions}\textup{(1)}. Hence the homology in degree $i$ is $D\otimes_AF_i$. Using $D\otimes_AA(-j)\cong D(-j)$, we obtain
$$
\Tor_i^A(D,M)\cong\bigoplus_jD(-j)^{b_{i,j}}\quad\text{and}\quad
b_{i,j}=\dim_D\bigl(\Tor_i^A(D,M)_j\bigr).
$$
\begin{definition}[Betti numbers]\label{def:betti-numbers}
We write $\beta^A_{i,j}(M)=b_{i,j}$ for the \emph{graded Betti numbers}, $\beta^A_i(M)=\sum_j\beta^A_{i,j}(M)$ for the total Betti numbers, and $\beta^A(M)$ for the Betti table; see also \cite[Definition~12.1]{Peeva}.
\end{definition}

Proposition~\ref{prop:syzygy-theorem} below shows that the minimal resolution is finite.
\begin{definition}\label{def:pd-reg}
For $M\ne0$, put
$$
\pd_A(M)=\max\{i\mid\beta^A_{i,j}(M)\ne0\text{ for some }j\}\quad\text{and}\quad \reg_A(M)=\max\{j-i\mid\beta^A_{i,j}(M)\ne0\}.
$$
The second invariant is the regularity defined from the minimal resolution, as in \cite[\S20.5]{Eisenbud}.
\end{definition}
Let $F_\bullet$ be the finite minimal graded free resolution of $M\ne0$, and set $s=\max\{i\mid F_i\ne0\}$. Since $F_s\ne0$, at least one $\beta^A_{s,j}(M)$ is positive. By the Tor formula,
$$
\Tor_s^A(D,M)\cong D\otimes_A F_s
\cong\bigoplus_j D(-j)^{\beta^A_{s,j}(M)}\ne0.
$$
A projective resolution of length less than $s$ would have zero homology in degree $s$ after tensoring with $D$. Since $F_\bullet$ has length $s$, we have $\pd_A(M)=s$, as in the commutative local case \cite[Corollary~19.5]{Eisenbud}. In terms of the Betti numbers, $\pd_A(M)=\max\{i\mid \beta^A_{i,j}(M)\ne0\text{ for some }j\}$.

If $\partial_{i+1}=0$ for some $i\ge0$, then exactness and minimality imply
$$
F_{i+1}=\ker\partial_{i+1}=\operatorname{im}\partial_{i+2}
\subseteq A_+F_{i+1}.
$$
By the graded Nakayama lemma, $F_{i+1}=0$. In particular, if $F_i=0$, then $\partial_{i+1}=0$, and induction shows that $F_j=0$ for all $j>i$.

\begin{proposition}\label{prop:syzygy-theorem}
Every finitely generated graded left $A$-module has a minimal graded free resolution of length at most $n$. In particular, $\pd_A(M)\le n$ for $M\ne0$.
\end{proposition}
\begin{proof}
The \emph{left global dimension} of a ring $R$ is $\operatorname{l.gl.dim}(R)=\sup_N\pd_R(N)$, where $N$ ranges over all left $R$-modules \cite[Remark after Theorem~4.1.2, p.~92]{Weibel}. Since every left $D$-module is free, $\operatorname{l.gl.dim}(D)=0$. By \cite[Theorem~4.3.7]{Weibel}, $\operatorname{l.gl.dim}(A)=n$ for $A=D[x_1,\dots,x_n]$ with central indeterminates. Hence $\Tor_i^A(D,M)=0$ for $i>n$. Let $F_\bullet\to M$ be a minimal graded free resolution, which exists by Proposition~\ref{prop:minimal-resolutions}\textup{(3)}. For $i>n$, the Tor formula implies
$$
0=\Tor_i^A(D,M)\cong D\otimes_A F_i
\cong\bigoplus_j D(-j)^{\beta^A_{i,j}(M)}.
$$
Each summand $A(-j)$ of $F_i$ contributes a nonzero summand $D(-j)$ to this direct sum. Thus $F_i=0$ for $i>n$. The resolution has length at most $n$, so $\pd_A(M)\le n$ when $M\ne0$.
\end{proof}

\begin{definition}\label{def:pure}
A minimal graded free resolution of length $s$ is \emph{pure of type} $(d_0,\dots,d_s)$ if $F_i\cong A(-d_i)^{\beta_i}$ with $\beta_i>0$ for $0\le i\le s$; see \cite[Definition~5.4, p.~426]{FloystadMcCulloughPeeva}. Then $d_0<\cdots<d_s$.

For $0\le i<s$, we have $F_{i+1}\ne0$. Since $\partial_{i+1}=0$ would imply $F_{i+1}=0$, the matrix of $\partial_{i+1}$ has a nonzero entry $a$. By Lemma~\ref{lem:matrix-degree}, $\deg a=d_{i+1}-d_i$. Since $a\in A_+$ by minimality, $d_{i+1}-d_i>0$. The Betti table determines $s$ and, for each $0\le i\le s$, the degrees $j$ with $\beta^A_{i,j}(M)\ne0$. Purity means that there is exactly one such degree in each of these homological degrees, and these degrees, in order, form the type.
\end{definition}

\begin{proposition}\label{prop:finite-length}
A finitely generated graded left $A$-module $M$ has finite length if and only if $\dim_DM<\infty$. Under these conditions, $H_M$ is a Laurent polynomial and $\dim_DM=H_M(1)$.
\end{proposition}
\begin{proof}
Suppose $\dim_D M<\infty$. For $A$-submodules $L\subsetneq N$ of $M$, we have $\dim_D L<\dim_D N$. Thus a strict chain from $0$ to $M$ has at most $\dim_D M$ inclusions. Choose a chain $0=U_0\subsetneq\cdots\subsetneq U_\ell=M$ of greatest length. A nonzero proper submodule of $U_j/U_{j-1}$ would lift to a submodule strictly between $U_{j-1}$ and $U_j$. Hence this chain is a composition series.

Conversely, suppose $M$ has finite length. The descending chain of submodules $A_+^rM$, starting with $A_+^0M=M$, stabilizes. Choose $r\ge0$ with $A_+^rM=A_+^{r+1}M$ and put $N=A_+^rM$. Then
$$
N=A_+^rM=A_+^{r+1}M=A_+N.
$$
The submodule $N$ is graded and finitely generated, since $M$ is noetherian. Proposition~\ref{prop:generators}\textup{(2)} therefore implies $N=0$. If $|\alpha|\ge r$, then $x^\alpha$ is divisible by a product of $r$ variables, so $x^\alpha\in A_+^r$ and $x^\alpha M=0$. If $m_1,\dots,m_k$ are homogeneous generators of $M$, the elements $x^\alpha m_i$ with $|\alpha|<r$ therefore span $M$ over $D$, so $M=\sum_{i=1}^k\sum_{|\alpha|<r}D x^\alpha m_i$. This is a finite sum, so $\dim_D M<\infty$.

Only finitely many of the components in $M=\bigoplus_j M_j$ are nonzero. Hence $H_M(t)$ is a Laurent polynomial, and $H_M(1)=\sum_j\dim_D M_j=\dim_D M$.
\end{proof}

\subsection{Changing generators}\label{subsec:changing-generators}
To compute the syzygies of another generating list for the same ideal, we express each list in terms of the other. The formulas below hold over any ring with identity. The corresponding formulas for left PBW rings and for commutative polynomial rings appear in \cite[\S4.2, (4.3)]{GomezTorrecillas} and \cite[Chapter~5, \S3, Proposition~(3.8)]{CLO}, respectively.

\begin{definition}[Syzygy]
Let $R$ be a ring with identity and $H\in R^{k\times1}$ a column. A \emph{left syzygy} of $H$ is a row $u\in R^{1\times k}$ with $uH=0$. These rows form the left $R$-module
$$
\Syz_{\ell,R}(H)=\{u\in R^{1\times k}\mid uH=0\}.
$$
For $R=A$, we write $\Syz_\ell(H)$, or $\Syz_\ell(h_1,\dots,h_k)$ when the column is understood. If $P$ is a matrix with $k$ columns, $\Row_R(P)\subseteq R^{1\times k}$ denotes the left submodule generated by its rows.
\end{definition}

\begin{proposition}\label{prop:transfer}
Let $R$ be a ring with identity, $F\in R^{s\times1}$, and $G\in R^{m\times1}$. Suppose $F=UG$ and $G=VF$, where $U\in M_{s\times m}(R)$ and $V\in M_{m\times s}(R)$. Then
$$
\begin{aligned}
\Syz_{\ell,R}(F)&=\Syz_{\ell,R}(G)V+\Row_R(I_s-UV),\\
\Syz_{\ell,R}(G)&=\Syz_{\ell,R}(F)U+\Row_R(I_m-VU).
\end{aligned}
$$
In particular, if $UV=I_s$, then $\Syz_{\ell,R}(F)=\Syz_{\ell,R}(G)V$. If $w_1,\dots,w_k$ generate $\Syz_{\ell,R}(G)$, then $w_1V,\dots,w_kV$ together with the rows of $I_s-UV$ generate $\Syz_{\ell,R}(F)$.
\end{proposition}
\begin{proof}
For $w\in\Syz_{\ell,R}(G)$, we have $(wV)F=wG=0$. Also $(I_s-UV)F=F-UG=0$, so every row of $I_s-UV$ belongs to $\Syz_{\ell,R}(F)$. Thus $\Syz_{\ell,R}(G)V+\Row_R(I_s-UV)\subseteq\Syz_{\ell,R}(F)$. Conversely, take $u\in\Syz_{\ell,R}(F)$. Then $(uU)G=uF=0$, and
$$
u=(uU)V+u(I_s-UV).
$$
The first summand belongs to $\Syz_{\ell,R}(G)V$. The second is a left $R$-linear combination of the rows of $I_s-UV$. This proves the reverse inclusion. The second equality follows by interchanging $F$ with $G$ and $U$ with $V$.

If $w_1,\dots,w_k$ generate $\Syz_{\ell,R}(G)$, then $\Syz_{\ell,R}(G)V=\sum_{j=1}^k R(w_jV)$ by left $R$-linearity of right multiplication by $V$. Thus the rows $w_jV$ and the rows of $I_s-UV$ generate $\Syz_{\ell,R}(F)$. When $UV=I_s$, the row module $\Row_R(I_s-UV)$ is zero, so $\Syz_{\ell,R}(F)=\Syz_{\ell,R}(G)V$.
\end{proof}

\subsection{Syzygies of a Gr\"obner basis}\label{subsec:schreyer}
We compute the syzygies of a left Gr\"obner basis. The formulas in \S~\ref{subsec:changing-generators} then determine the syzygies of another generating list of the same ideal.

\begin{definition}[Gr\"obner basis]
Let $J\subseteq A$ be a left ideal. A finite ordered list $H=(h_1,\dots,h_m)$ of nonzero elements of $J$, with repetitions allowed, is a \emph{left Gr\"obner basis} of $J$ if, for every nonzero $r\in J$, some $\LM(h_i)$ divides $\LM(r)$. The \emph{initial ideal} of $J$ with respect to $\prec$ is $\operatorname{in}_\prec(J)=\sum_{0\ne r\in J}A\LM(r)$.
\end{definition}

Let $G=(g_1,\dots,g_m)^t$ be a column of monic polynomials and let $e_1,\dots,e_m$ be the canonical basis of $A^{1\times m}$. A \emph{module monomial} is an element $x^\alpha e_i$ which divides $x^\beta e_j$ if $i=j$ and $x^\alpha\mid x^\beta$. As in the commutative case \cite[Chapter~5, \S3, Theorem~(3.3)]{CLO}, define the \emph{Schreyer order} by setting $x^\alpha e_i\prec_Gx^\beta e_j$ if $x^\alpha\LM(g_i)\prec x^\beta\LM(g_j)$, or if these products are equal and $i>j$. Write $\LM_G$, $\LC_G$, and $\LT_G$ for the leading monomial, leading coefficient, and leading term in this order. The coefficients are written on the left. A left Gr\"obner basis of a submodule is defined by divisibility of leading module monomials.

\begin{lemma}\label{lem:schreyer-order}
The order $\prec_G$ is a monomial order on module monomials.
\end{lemma}
\begin{proof}
For a module monomial $u=x^\alpha e_i$, put $w(u)=x^\alpha\LM(g_i)$. The pair $(w(u),-i)$ determines $u$, since $x^\alpha=w(u)/\LM(g_i)$. Lexicographic comparison of these pairs defines the total order $\prec_G$. Multiplication by $x^\gamma$ fixes the basis index and satisfies $w(x^\gamma u)=x^\gamma w(u)$. By multiplicativity of $\prec$, $u\prec_G v$ implies $x^\gamma u\prec_G x^\gamma v$.

Let $U$ be a nonempty set of module monomials. Choose $u_0=x^{\alpha_0}e_{i_0}\in U$ with $w(u_0)=\min_{\prec}\{w(u)\mid u\in U\}$, taking the largest possible $i_0$. These choices are possible since $\prec$ is a well-order and the indices lie in $\{1,\dots,m\}$. For $u=x^\alpha e_i\in U$, we have $w(u_0)\preceq w(u)$, with $i_0\ge i$ if $w(u_0)=w(u)$. Thus $u_0\preceq_G u$ for every $u\in U$.
\end{proof}

For $0\ne a\in A$ and $0\ne s\in A^{1\times m}$, let $\LT(a)=cx^\gamma$ and $\LT_G(s)=dx^\alpha e_i$. Each product of monomials $x^\eta$ of $a$ and $x^\beta e_j$ of $s$ satisfies $x^{\eta+\beta}e_j\preceq_G x^{\eta+\alpha}e_i\preceq_G x^{\gamma+\alpha}e_i$. Equality holds throughout exactly when $\eta=\gamma$ and $(\beta,j)=(\alpha,i)$. The coefficient of $x^{\gamma+\alpha}e_i$ in $as$ is therefore $cd\ne0$, so $\LT_G(as)=\LT(a)\LT_G(s)$.

\begin{lemma}\label{lem:left-division}
Let $g_1,\dots,g_m\in A$ be monic, and let $f\in A$. There is an expression $f=\sum_kq_kg_k+\rho$ such that no monomial of $\rho$ is divisible by any $\LM(g_k)$. If $f\ne0$, then $\LM(q_k)\LM(g_k)\preceq\LM(f)$ whenever $q_k\ne0$. For $f=0$, take all $q_k$ and $\rho$ to be zero. If the $g_k$ form a left Gr\"obner basis of $J$ and $f\in J$, then $\rho=0$. In particular, a left Gr\"obner basis of $J$ generates $J$.
\end{lemma}
\begin{proof}
For $f=0$, take all $q_k$ and $\rho$ to be zero. Suppose $f\ne0$. Start with $h=f$ and $q_1=\cdots=q_m=\rho=0$. The left division algorithm \cite[\S4.1, Algorithm~3]{GomezTorrecillas} implies the equality
$$
f=\sum_{k=1}^m q_k g_k+h+\rho.
$$
If $h\ne0$, write $\LT(h)=cx^\gamma$. When $\LM(g_k)=x^{\delta_k}$ divides $x^\gamma$, add $cx^{\gamma-\delta_k}$ to $q_k$ and subtract $cx^{\gamma-\delta_k}g_k$ from $h$. Since $g_k$ is monic, the subtracted polynomial has leading term $cx^\gamma$. If no $\LM(g_k)$ divides $x^\gamma$, subtract $cx^\gamma$ from $h$ and add it to $\rho$.

Each step preserves the equality and strictly decreases $\LM(h)$ whenever the new $h$ is nonzero. Since $\prec$ is a well-order, the process ends with $h=0$. Every monomial of $\rho$ comes from a step in which no $\LM(g_k)$ divides it. Each monomial $x^\eta$ of the final $q_k$ was added at a step with $x^\eta\LM(g_k)=\LM(h)\preceq\LM(f)$ before the subtraction. In particular, $\LM(q_k)\LM(g_k)\preceq\LM(f)$ whenever $q_k\ne0$.

If the $g_k$ form a left Gr\"obner basis of $J$ and $f\in J$, then $\rho=f-\sum_k q_k g_k\in J$. For $\rho\ne0$, some $\LM(g_k)$ would divide $\LM(\rho)$ by the Gr\"obner basis property. The construction of $\rho$ excludes this, so $\rho=0$. Hence $J=\sum_k A g_k$.
\end{proof}

The \emph{least common multiple} of $x^\alpha$ and $x^\beta$ is $\operatorname{lcm}(x^\alpha,x^\beta)=x^\gamma$, where $\gamma_i=\max(\alpha_i,\beta_i)$. It divides every common monomial multiple. For $i<j$, put $m_{ij}=\operatorname{lcm}(\LM(g_i),\LM(g_j))$, and define $\mu_{ij},\nu_{ij}$ by $\mu_{ij}\LM(g_i)=m_{ij}=\nu_{ij}\LM(g_j)$. The \emph{left $S$-polynomial} is $S_\ell(g_i,g_j)=\mu_{ij}g_i-\nu_{ij}g_j$. As in the commutative definition \cite[Chapter~2, \S9, Definition~5, p.~107]{CLOIVA}, an \emph{lcm representation} of $S_\ell(g_i,g_j)$ is an expression
$$
S_\ell(g_i,g_j)=\sum_{k=1}^mq_{ijk}g_k,\qquad \LM(q_{ijk})\LM(g_k)\prec m_{ij}\quad\text{whenever }q_{ijk}\ne0.
$$

When the $g_k$ form a left Gr\"obner basis, Lemma~\ref{lem:left-division} supplies an lcm representation, since $S_\ell(g_i,g_j)=0$ or $\LM(S_\ell(g_i,g_j))\prec m_{ij}$. Fix one such representation for each pair and put $s_{ij}=\mu_{ij}e_i-\nu_{ij}e_j-\sum_kq_{ijk}e_k$.

\begin{proposition}\label{prop:schreyer}
Let $G$ be a monic left Gr\"obner basis of the left ideal generated by its entries, with an lcm representation fixed for each pair $i<j$. Then $\{s_{ij}\}_{i<j}$ is a monic left Gr\"obner basis of $\Syz_\ell(G)$ with respect to $\prec_G$, and $\Syz_\ell(G)=\sum_{i<j}As_{ij}$. If $m=1$, this module is zero.
\end{proposition}
\begin{proof}
The construction is the one in \cite[\S4.2, Corollary~6]{GomezTorrecillas}. Here $s_{ij}G=0$ and the bounds in the chosen lcm representation put every term of $\sum_kq_{ijk}e_k$ strictly below $\mu_{ij}e_i$ in $\prec_G$. Since $\mu_{ij}\LM(g_i)=\nu_{ij}\LM(g_j)=m_{ij}$ and $i<j$, we have $\LT_G(s_{ij})=\mu_{ij}e_i$.

Let $0\ne s=(a_1,\ldots,a_m)\in\Syz_\ell(G)$ and put $T=\max_{a_k\ne0}\LM(a_k)\LM(g_k)$. The coefficient of $T$ in $sG$ is zero. Since each summand whose leading monomial is $T$ contributes a nonzero coefficient, at least two indices attain $T$. Let $i$ be the least of them and choose another index $j>i$ attaining $T$. Both $\LM(g_i)$ and $\LM(g_j)$ divide $T$, so $m_{ij}\mid T$. By the tie rule,
$$
\LM_G(s)=\LM(a_i)e_i=(T/m_{ij})\mu_{ij}e_i=(T/m_{ij})\LM_G(s_{ij}).
$$
Subtracting $\LC(a_i)(T/m_{ij})s_{ij}$ cancels $\LT_G(s)$ and leaves a syzygy. The well-order from Lemma~\ref{lem:schreyer-order} ensures that repeated reduction ends at zero. Thus the $s_{ij}$ generate $\Syz_\ell(G)$ and their leading module monomials satisfy the required divisibility. If $m=1$, there can be no nonzero syzygy since $A$ is a domain and $g_1\ne0$.
\end{proof}

\begin{proposition}\label{prop:strict-criterion}
Let $H=(h_1,\dots,h_m)^t$ be a column of monic polynomials and put $J=\sum_iAh_i$. If every $S_\ell(h_i,h_j)$ with $i<j$ has an lcm representation in terms of $H$, then $H$ is a monic left Gr\"obner basis of $J$.
\end{proposition}
\begin{proof}
The commutative version of Buchberger's criterion uses lcm representations \cite[Chapter~2, \S9, Theorem~6]{CLOIVA}. For left PBW rings, the criterion is stated in terms of reduction of the $S$-polynomials to zero \cite[\S4.2, Theorem~7]{GomezTorrecillas}. For $i<j$, put $m_{ij}=\operatorname{lcm}(\LM(h_i),\LM(h_j))$. For $0\ne r\in J$, choose an expression $r=\sum_k a_kh_k$ for which $T=\max_{a_k\ne0}\LM(a_k)\LM(h_k)$ is least. Let $\Lambda$ be the set of indices attaining $T$. For $k\in\Lambda$ write $\LT(a_k)=c_kn_k$. If $\sum_{k\in\Lambda}c_k\ne0$, then $\LM(r)=T$ and $\LM(h_k)\mid\LM(r)$ for $k\in\Lambda$.

Suppose $\sum_{k\in\Lambda}c_k=0$ and let $i=\min\Lambda$. Then $m_{ij}\mid T$ for $j\in\Lambda\smallsetminus\{i\}$, and
$$
\sum_{k\in\Lambda}c_kn_kh_k
=-\sum_{j\in\Lambda\smallsetminus\{i\}}c_j(T/m_{ij})S_\ell(h_i,h_j).
$$
Write $S_\ell(h_i,h_j)=\sum_kq_{ijk}h_k$ using the chosen lcm representation. Then $\LM(q_{ijk})\LM(h_k)\prec m_{ij}$ whenever $q_{ijk}\ne0$. Thus $\LM(c_j(T/m_{ij})q_{ijk})\LM(h_k)\prec T$. Subtract $c_kn_k$ from $a_k$ for $k\in\Lambda$ and substitute the displayed sum. The resulting expression $r=\sum_k a'_kh_k$ has $\max_{a'_k\ne0}\LM(a'_k)\LM(h_k)\prec T$, contrary to the choice of $T$. Thus the nonzero coefficient case applies.
\end{proof}

\begin{proposition}\label{prop:standard-monomials}
The monomials outside $\operatorname{in}_\prec(J)$ are the \emph{standard monomials} of $J$ with respect to $\prec$. If $J$ has a monic left Gr\"obner basis, their classes form a left $D$-basis of $A/J$.
\end{proposition}
\begin{proof}
Let $H=(h_1,\dots,h_m)^t$ be a monic left Gr\"obner basis of $J$. Since $h_i\in J$, each $A\LM(h_i)$ is contained in $\operatorname{in}_\prec(J)$. Conversely, for every $0\ne f\in J$, some $\LM(h_i)$ divides $\LM(f)$, so $A\LM(f)\subseteq A\LM(h_i)$. Thus $\operatorname{in}_\prec(J)=\sum_{i=1}^m A\LM(h_i)$. By the monomial basis of $A$, the monomials in this ideal are exactly those divisible by at least one $\LM(h_i)$. Let $W$ be the left $D$-span of the standard monomials.

By Lemma~\ref{lem:left-division}, each $f\in A$ can be written as $f=\sum_i q_i h_i+\rho$, with $\rho\in W$. Hence $A=J+W$. If $0\ne r\in J\cap W$, then $\LM(r)\in\operatorname{in}_\prec(J)$ because $r\in J$. But $\LM(r)$ is a standard monomial because $r\in W$. Thus $J\cap W=0$, and the quotient map restricts to an isomorphism of left $D$-vector spaces
$$
W\xrightarrow{\ \sim\ }A/J\quad\text{defined by}\quad w\longmapsto w+J.
$$
The standard monomials form a basis of $W$, so their classes form a left $D$-basis of $A/J$, as in the case of a polynomial ring over a field \cite[Chapter~2, \S2]{CLO}.
\end{proof}

\section{Syzygies and minimal resolutions}\label{sec:dichotomy}\label{pg:s3}
Let $D$ be a division ring and let $A=D[x,y,z]$, where $x,y,z$ are central indeterminates of degree $1$. Fix a monomial order with $z\prec y\prec x$. For positive integers $p,q$ and $a,b\in D$, put
$$
f=x^p+az^p,\qquad g=y^q+bz^q,\qquad
\kappa=ba-ab,\qquad I=Af+Ag,\qquad M=A/I.
$$
If $\kappa\ne0$, the factor $\kappa^{-1}$ need not be central, so the order of the factors in products such as $f\kappa^{-1}g$ matters.

Since $f$ and $g$ have positive degrees, $I\subseteq A_+$ and $M/A_+M\cong A/A_+\cong D$ in degree zero. Hence $\beta^A_{0,0}(M)=1$. The order $z\prec y\prec x$ implies $z^p\prec x^p$ and $z^q\prec y^q$. Thus $f$ and $g$ are monic with leading monomials $x^p$ and $y^q$. In the expansion of $gf-fg$, the terms involving $x^p$ or $y^q$ cancel, so $gf-fg=\kappa z^{p+q}$. Since $x^p$ and $y^q$ are relatively prime, their least common multiple is $x^py^q$, and
$$
S_\ell(f,g)=y^qf-x^pg
=gf-fg-(bz^q)f+(az^p)g
=(-bz^q)f+(az^p)g+\kappa z^{p+q}.
$$

\begin{proposition}\label{prop:commuting-case}
If $\kappa=0$, then $\{f,g\}$ is a monic left Gr\"obner basis of $I$ and $\operatorname{in}_\prec(I)=Ax^p+Ay^q$.
\end{proposition}
\begin{proof}
When $\kappa=0$, the formula for the left $S$-polynomial reduces to $S_\ell(f,g)=(-bz^q)f+(az^p)g$. By multiplying $z^q\prec y^q$ by $x^p$ and $z^p\prec x^p$ by $y^q$, we obtain
$$
\begin{aligned}
\LM(-bz^q)\LM(f)&=x^pz^q\prec x^py^q &&\text{if } b\ne0,\\
\LM(az^p)\LM(g)&=z^py^q\prec x^py^q &&\text{if } a\ne0.
\end{aligned}
$$
These inequalities show that the expression for $S_\ell(f,g)$ is an lcm representation. By Proposition~\ref{prop:strict-criterion}, $\{f,g\}$ is a monic left Gr\"obner basis of $I$. Hence $\operatorname{in}_\prec(I)=A\LM(f)+A\LM(g)=Ax^p+Ay^q$.
\end{proof}

\begin{proposition}\label{prop:noncommuting-case}
If $\kappa\ne0$, put $h=z^{p+q}$. Then $h=(\kappa^{-1}g)f-(\kappa^{-1}f)g\in I$, $f,g,h$ form a monic left Gr\"obner basis of $I$, and $\operatorname{in}_\prec(I)=Ax^p+Ay^q+Ah$.
\end{proposition}
\begin{proof}
Suppose $\kappa\ne0$ and put $h=z^{p+q}$. It follows from $gf-fg=\kappa h$ that $h=(\kappa^{-1}g)f-(\kappa^{-1}f)g$. Both coefficients have positive degree, so $h\in A_+I\subseteq I$. Since $\LM(f)\nmid h$ and $\LM(g)\nmid h$, $\{f,g\}$ is not a left Gr\"obner basis of $I$.

The polynomials $f,g,h$ generate $I$. Their leading monomials $x^p,y^q,h$ are pairwise relatively prime, so the least common multiples for the three pairs are $x^py^q$, $x^ph$, and $y^qh$. The corresponding left $S$-polynomials are
$$
\begin{aligned}
S_\ell(f,g)&=(-bz^q)f+(az^p)g+\kappa h,\\
S_\ell(f,h)&=hf-x^ph=(az^p)h,\\
S_\ell(g,h)&=hg-y^qh=(bz^q)h.
\end{aligned}
$$

Since $\kappa\ne0$, both $a$ and $b$ are nonzero. The leading monomials of the three summands in the expression for $S_\ell(f,g)$ are $x^pz^q$, $z^py^q$, and $h$. By multiplicativity of the monomial order, $z^p\prec x^p$ and $z^q\prec y^q$ imply
$$
h=z^pz^q\prec x^pz^q\prec x^py^q
\quad\text{and}\quad
z^py^q\prec x^py^q.
$$
The leading monomial of $(az^p)h$ is $z^ph\prec x^ph$, and that of $(bz^q)h$ is $z^qh\prec y^qh$. Thus all three expressions are lcm representations. By Proposition~\ref{prop:strict-criterion}, $\{f,g,h\}$ is a monic left Gr\"obner basis of $I$. Its leading monomials generate $\operatorname{in}_\prec(I)=Ax^p+Ay^q+Ah$.
\end{proof}

\Needspace{4\baselineskip}
\begin{theorem}\label{thm:dichotomy}
With the notation above, the following hold.
\begin{enumerate}[label=\textup{(\arabic*)}]
\item If $\kappa=0$, then $\{f,g\}$ is a monic left Gr\"obner basis of $I$, and $M$ has the minimal graded free resolution
$$
0\longrightarrow A(-(p+q))\xrightarrow{\ \cdot(g,-f)\ }A(-p)\oplus A(-q)
\xrightarrow{\ \cdot(f,g)^t\ }A\longrightarrow M\longrightarrow0,
$$
with $H_M(t)=(1-t^p)(1-t^q)(1-t)^{-3}$, $\pd_A(M)=2$, $\reg_A(M)=p+q-2$, and $\dim_DM=\infty$. If $p\ne q$, the nonzero graded Betti numbers are $\beta^A_{0,0}=1$, $\beta^A_{1,p}=\beta^A_{1,q}=1$, and $\beta^A_{2,p+q}=1$. If $p=q=d$, they are $\beta^A_{0,0}=1$, $\beta^A_{1,d}=2$, and $\beta^A_{2,2d}=1$. The total Betti numbers are $(1,2,1)$.
\item If $\kappa\ne0$, put $h=z^{p+q}$. The polynomials $f,g,h$ form a monic left Gr\"obner basis of $I$, and $M$ has the minimal graded free resolution
$$
0\to A(-2p-2q)\xrightarrow{\cdot P_3} A(-(2p+q))\oplus A(-(p+2q))
\xrightarrow{\cdot P_2} A(-p)\oplus A(-q)\xrightarrow{\cdot P_1} A\to M\to0,
$$
where
$$
P_1=\begin{pmatrix}f\\g\end{pmatrix},\qquad
P_3=\begin{pmatrix}-\kappa^{-1}g & \kappa^{-1}f\end{pmatrix},\qquad
P_2=\begin{pmatrix}
h-f\kappa^{-1}g & f\kappa^{-1}f\\[2pt]
-g\kappa^{-1}g & h+g\kappa^{-1}f
\end{pmatrix}.
$$
Moreover, $H_M(t)=(1-t^p)(1-t^q)(1-t^{p+q})(1-t)^{-3}$, $\dim_DM=pq(p+q)$, $\pd_A(M)=3$, and $\reg_A(M)=2p+2q-3$. If $p\ne q$, the nonzero graded Betti numbers are $\beta^A_{0,0}=1$, $\beta^A_{1,p}=\beta^A_{1,q}=1$, $\beta^A_{2,2p+q}=\beta^A_{2,p+2q}=1$, and $\beta^A_{3,2p+2q}=1$. If $p=q=d$, they are $\beta^A_{0,0}=1$, $\beta^A_{1,d}=2$, $\beta^A_{2,3d}=2$, and $\beta^A_{3,4d}=1$. The total Betti numbers are $(1,2,2,1)$.
\end{enumerate}
\end{theorem}

\begin{proof}
\emph{The case $\kappa=0$.} Put $C=(f,g)^t$. By Proposition~\ref{prop:commuting-case}, $\{f,g\}$ is a monic left Gr\"obner basis of $I$. In view of the lcm representation $S_\ell(f,g)=(-bz^q)f+(az^p)g$, Proposition~\ref{prop:schreyer} shows that $\Syz_\ell(C)$ is generated by $s_{12}=y^qe_1-x^pe_2-\bigl((-bz^q)e_1+(az^p)e_2\bigr)=(g,-f)$. The entries of $C$ have degrees $p,q$, and those of $s_{12}$ have degrees $q=(p+q)-p$ and $p=(p+q)-q$. Hence both maps have degree zero by Lemma~\ref{lem:matrix-degree}.

Since $s_{12}C=gf-fg=0$, the maps form a complex. Exactness at $A(-p)\oplus A(-q)$ and at $A$ follows from $\operatorname{im}\rho_{s_{12}}=As_{12}=\ker\rho_C$ and $\operatorname{im}\rho_C=I=\ker(A\to M)$. If $us_{12}=0$, then $ug=0$, so $u=0$ by Corollary~\ref{cor:elementary}(1) and $g\ne0$. Together with surjectivity of $A\to M$, this proves exactness. The entries of $C$ and $s_{12}$ lie in $A_+$ since $p,q\ge1$, so the resolution is minimal by Proposition~\ref{prop:minimal-resolutions}(1).

By Proposition~\ref{prop:standard-monomials}, the classes of $x^iy^jz^k$ for $0\le i<p$, $0\le j<q$, and $k\ge0$ form a left $D$-basis of $M$. This basis is infinite. Counting its elements by degree, we obtain
$$
H_M(t)=
\left(\sum_{i=0}^{p-1}t^i\right)
\left(\sum_{j=0}^{q-1}t^j\right)
\left(\sum_{k\ge0}t^k\right)
=\frac{(1-t^p)(1-t^q)}{(1-t)^3}.
$$
By Definition~\ref{def:betti-numbers}, the graded Betti numbers count the shifts in the lists $(0)$, $(p,q)$, and $(p+q)$ in their respective homological degrees. The two shifts in homological degree one contribute separately if $p\ne q$ and with multiplicity two at $d$ if $p=q=d$. The total Betti numbers are $(1,2,1)$, so $\pd_A(M)=2$. Since $p,q\ge1$, we have $\reg_A(M)=\max\{0,p-1,q-1,p+q-2\}=p+q-2$.

\emph{The case $\kappa\ne0$.} Put $h=z^{p+q}$ and $G=(f,g,h)^t$. By Proposition~\ref{prop:noncommuting-case}, the entries of $G$ form a monic left Gr\"obner basis of $I$. Write $C=P_1=(f,g)^t$ and $R=P_3=(-\kappa^{-1}g,\kappa^{-1}f)$.

\emph{Step 1. The complex and its grading.} The matrices satisfy $P_2=hI_2+CR$ and $RC=\kappa^{-1}(fg-gf)=-h$. Since $h$ is central,
$$
\begin{aligned}
P_2P_1&=hC+C(RC)=hC-Ch=0,\\
P_3P_2&=Rh+(RC)R=Rh-hR=0.
\end{aligned}
$$
The entries of $P_1$ lie in $A_p,A_q$, in that order, and those of $P_3$ lie in $A_q,A_p$. Since $h\in A_{p+q}$ and $\kappa^{-1}\in A_0$, the entries of $P_2$ lie in the homogeneous components with indices
$$
\begin{pmatrix}
p+q&2p\\
2q&p+q
\end{pmatrix}
=
\begin{pmatrix}
(2p+q)-p &(2p+q)-q\\
(p+2q)-p &(p+2q)-q
\end{pmatrix}.
$$
These indices and those for $P_1,P_3$ are the source shifts minus the target shifts. Thus all three maps have degree zero by Lemma~\ref{lem:matrix-degree}.

\emph{Step 2. Exactness at $F_1=A(-p)\oplus A(-q)$.} To compute $\ker\rho_{P_1}$, take
$$
U=\begin{pmatrix}1&0&0\\0&1&0\end{pmatrix},
\qquad
V=\begin{pmatrix}1&0\\0&1\\\kappa^{-1}g&-\kappa^{-1}f\end{pmatrix}.
$$
Then $C=UG$ and $UV=I_2$. The last entry of $VC$ is $(\kappa^{-1}g)f-(\kappa^{-1}f)g=h$, so $G=VC$. The entries in the last row have degrees $q=(p+q)-p$ and $p=(p+q)-q$. Hence $\rho_V\colon F_1\oplus A(-(p+q))\to F_1$ has degree zero by Lemma~\ref{lem:matrix-degree}.

By Proposition~\ref{prop:schreyer}, the syzygies associated to the lcm representations in Proposition~\ref{prop:noncommuting-case} generate $\Syz_\ell(G)$. They are
$$
s_{12}=(g,-f,-\kappa),\qquad
s_{13}=(h,0,-f),\qquad
s_{23}=(0,h,-g).
$$ Their images under right multiplication by $V$ are
$$
\begin{aligned}
s_{12}V&=(0,0),\\
s_{13}V&=(h-f\kappa^{-1}g,f\kappa^{-1}f),\\
s_{23}V&=(-g\kappa^{-1}g,h+g\kappa^{-1}f).
\end{aligned}
$$
The last two images are the rows of $P_2$. It follows from $UV=I_2$ and Proposition~\ref{prop:transfer} that
$$
\ker\rho_{P_1}=\Syz_\ell(C)=\Syz_\ell(G)V
=\Row_A(P_2)=\operatorname{im}\rho_{P_2}.
$$

\emph{Step 3. Exactness at $F_2$.} The identity $P_3P_2=0$ implies $AR\subseteq\ker\rho_{P_2}$. For the reverse inclusion, take $s=(u,v)\in\ker\rho_{P_2}$ and put $t=sC$. Then $0=sP_2=hs+tR$. In the first coordinate, this equality reads $hu=t\kappa^{-1}g$.

Since $h$ is central, $Ah=hA$ is a two-sided ideal, and $A/Ah\cong B[y]$ for $B=D[x,z]/(z^{p+q})$. Write a bar for reduction modulo $Ah$. The coordinate identity reduces to
$$
\overline t(\kappa^{-1}\overline g)=0,
\qquad
\kappa^{-1}\overline g=\kappa^{-1}y^q+\kappa^{-1}b\overline z^{\,q}.
$$
The leading coefficient of $\kappa^{-1}\overline g$ as a polynomial in $y$ is $\kappa^{-1}\in B^\times$. By Remark~\ref{lem:unit-leading-regular}(1), right multiplication by this polynomial is injective, so $\overline t=0$. Write $t=hw$ with $w\in A$. Then $0=hs+tR=h(s+wR)$. Since $h\ne0$ and $A$ is a domain, $s=-wR$. Thus $\ker\rho_{P_2}=AR=\operatorname{im}\rho_{P_3}$.

\emph{Step 4. Injectivity of $\rho_{P_3}$ and minimality.} If $uR=0$, then $u\kappa^{-1}g=0$ by the first coordinate, so $u=0$ since $\kappa^{-1}g\ne0$ and $A$ is a domain. At $A$, $\operatorname{im}\rho_{P_1}=Af+Ag=I=\ker(A\to M)$, and the quotient map is surjective. Thus the sequence is exact. All entries of $P_1,P_2,P_3$ lie in $A_+$ by the degree computations in Step~1 and $p,q\ge1$. Thus the resolution is minimal by Proposition~\ref{prop:minimal-resolutions}(1). By Proposition~\ref{prop:generators}(4), $f,g$ form a minimal homogeneous generating set of $I$.

\emph{Step 5. Betti numbers and Hilbert series.} By Proposition~\ref{prop:noncommuting-case}, the initial ideal is $Ax^p+Ay^q+Ah$. The standard monomials of $I$ are the $x^iy^jz^k$ with $0\le i<p$, $0\le j<q$, and $0\le k<p+q$. Their classes form a left $D$-basis of $M$ by Proposition~\ref{prop:standard-monomials}. Counting these $pq(p+q)$ classes by degree, we obtain
$$
H_M(t)=
\left(\sum_{i=0}^{p-1}t^i\right)
\left(\sum_{j=0}^{q-1}t^j\right)
\left(\sum_{k=0}^{p+q-1}t^k\right)
=\frac{(1-t^p)(1-t^q)(1-t^{p+q})}{(1-t)^3}.
$$
The shifts in homological degrees $0,1,2,3$ are $(0)$, $(p,q)$, $(2p+q,p+2q)$, and $(2p+2q)$, respectively. For $p\ne q$, the two shifts in each of degrees one and two are distinct. For $p=q=d$, they have multiplicity two at $d$ and $3d$, respectively. Thus the graded Betti numbers are as stated, with totals $(1,2,2,1)$ and $\pd_A(M)=3$. Finally, $p,q\ge1$ implies
$$
\begin{aligned}
\reg_A(M)
&=\max\{0,p-1,q-1,2p+q-2,p+2q-2,2p+2q-3\}\\
&=2p+2q-3.
\end{aligned}
$$
\end{proof}

For $\kappa\ne0$, we have $\operatorname{im}\rho_{P_2}=\Syz_\ell(f,g)$ by exactness. Since the resolution is minimal, Proposition~\ref{prop:generators}\textup{(4)} shows that the rows $s_{13}V$ and $s_{23}V$ of $P_2$ form a minimal homogeneous generating set of $\Syz_\ell(f,g)$.

The module in Theorem~\ref{thm:dichotomy}(2) has projective dimension three and its defining left ideal has two minimal generators. The following lemma describes the corresponding quotients over a field.

\begin{lemma}\label{lem:two-generated}
Let $S=K[x_1,\dots,x_N]$, where $K$ is a field and $N\ge1$, and let $M$ be a finitely generated graded $S$-module with $\beta^S_0(M)=\beta^S_{0,0}(M)=1$ and $\beta^S_1(M)=2$. Then $\pd_S(M)=2$ and $\beta^S_2(M)=1$.
\end{lemma}
\begin{proof}
Let $F_\bullet\to M$ be a minimal graded free resolution. The Betti numbers imply that $F_0\cong S$ and $F_1\cong S(-d)\oplus S(-d')$ for some integers $d,d'$. Choose homogeneous bases and write $\partial_1\colon S(-d)\oplus S(-d')\to S$ for the map defined by $\partial_1(r,s)=ru+sv$, where $u\in S_d$ and $v\in S_{d'}$. Exactness implies $M\cong S/(u,v)$, and $u,v$ minimally generate $(u,v)$ by Proposition~\ref{prop:generators}\textup{(4)}. In particular, $u,v\ne0$ and $d,d'>0$.

Let $c$ be a greatest common divisor of $u,v$ and put $e=\deg c$. Then $c$ is homogeneous, since factors of nonzero homogeneous polynomials in $S$ are homogeneous. Write $u=cu'$ and $v=cv'$, with $\gcd(u',v')=1$. Then $\deg u'=d-e$ and $\deg v'=d'-e$. If $e=d$, then $u'\in K^\times$ and $v=u(u')^{-1}v'\in(u)$, contradicting minimality. Likewise, $e=d'$ would imply $u\in(v)$. Thus $0\le e<\min\{d,d'\}$, so $u',v'\in S_+$.

For every $w\in S$, the condition $\gcd(u',v')=1$ implies $v'w\in(u')\Rightarrow w\in(u')$. Hence multiplication by $v'$ on $S/(u')$ is injective. Since $S$ is a domain and $(u',v')\subseteq S_+$, the pair $u',v'$ is a regular sequence. Its Koszul complex computes the syzygy module as $S(v',-u')$ \cite[Corollary~17.5]{Eisenbud}. Cancelling $c\ne0$ in $ru+sv=c(ru'+sv')$, we obtain
$$
\ker\partial_1=\{(r,s)\in S^2\mid ru'+sv'=0\}=S(v',-u').
$$

Put $\sigma=d+d'-e$. The map $\partial_2\colon S(-\sigma)\to S(-d)\oplus S(-d')$ defined by $\partial_2(t)=(tv',-tu')$ has degree $0$, since $\deg v'=\sigma-d$ and $\deg u'=\sigma-d'$. If $\partial_2(t)=0$, then $tv'=0$, so $t=0$. Together with the kernel calculation, this proves exactness of
$$
0\longrightarrow S(-\sigma)\xrightarrow{\partial_2}
S(-d)\oplus S(-d')\xrightarrow{\partial_1}S\longrightarrow M\longrightarrow0.
$$
The entries $u,v,v',-u'$ all lie in $S_+$, so this resolution is minimal by Proposition~\ref{prop:minimal-resolutions}\textup{(1)}. The term in homological degree $2$ is $S(-\sigma)$, so $\pd_S(M)=2$ and $\beta^S_2(M)=1$.
\end{proof}

\begin{corollary}\label{prop:notcomm}
Let $K$ be a field and let $N\ge1$. No Betti table in Theorem~\ref{thm:dichotomy}(2) is the Betti table of a finitely generated graded $K[x_1,\dots,x_N]$-module.
\end{corollary}
\begin{proof}
Put $S=K[x_1,\dots,x_N]$, and suppose a finitely generated graded $S$-module $N'$ has the Betti table in Theorem~\ref{thm:dichotomy}(2) for some $p,q$. In the minimal resolution $F_\bullet\to N'$, the table specifies $F_0\cong S$ and $F_1\cong S(-p)\oplus S(-q)$. Thus $\beta^S_0(N')=\beta^S_{0,0}(N')=1$ and $\beta^S_1(N')=2$. By Lemma~\ref{lem:two-generated}, $\pd_S(N')=2$. Hence $F_3=0$, contradicting the summand $S(-2p-2q)$ required by the entry $\beta^S_{3,2p+2q}(N')=1$.
\end{proof}

\section{The quaternion case}\label{sec:pure}\label{pg:s4}
Let $\mathbb H=\mathbb R\oplus\mathbb R i\oplus\mathbb R j\oplus\mathbb R k$ be the real quaternion division ring, with $i^2=j^2=k^2=ijk=-1$. Conjugation is defined by $\overline{a_0+a_1i+a_2j+a_3k}=a_0-a_1i-a_2j-a_3k$, where $a_0,\dots,a_3\in\mathbb R$.

Let $R=\mathbb R[t_{r\ell}\mid r=1,2,\ 0\le\ell\le3]$ and put $\mathcal A=\mathbb H\otimes_{\mathbb R}R$. For $r=1,2$, set $q_r=t_{r0}+it_{r1}+jt_{r2}+kt_{r3}\in\mathcal A$.

Let $U_r\in M_4(R)$ be the matrix of left multiplication by $\bar q_r$ on $\mathcal A$ in the $R$-basis $(1,i,j,k)$, and put $M_2=\operatorname{coker}(U_1\ U_2)$. Adams and Loustaunau compute a minimal graded free resolution of $\operatorname{im}(U_1\ U_2)$ \cite[equation~(5), p.~8]{AdamsLoustaunau}. Appending its inclusion into $R^4$ and the quotient map, we obtain the minimal graded free resolution
$$
0\longrightarrow R(-4)^4\longrightarrow R(-3)^8\longrightarrow
R(-1)^8\xrightarrow{(U_1\ U_2)}R^4\longrightarrow M_2\longrightarrow0,
$$
where $\deg t_{r\ell}=1$. The shifts and ranks also follow from \cite[equation~(9) and Corollary~3.6, p.~13]{AdamsLoustaunau}.

\begin{proposition}\label{rem:quaternion-specialization}
In the ordinary grading, the resolution of $M_2$ recalled above is isomorphic to the restriction to $R$ of a minimal graded free resolution of the left $\mathcal A$-module $\mathcal M=\mathcal A/(\mathcal A q_1+\mathcal A q_2)$. Let $p,q$ be positive integers and let $a,b\in\mathbb H$ satisfy $\kappa=ba-ab\ne0$. Assign degree $p$ to the coordinates of $q_1$ and degree $q$ to those of $q_2$. Let $\varphi\colon\mathcal A\to A=\mathbb H[x,y,z]$ be the homomorphism fixing $\mathbb H$ and defined by the substitutions
$$
q_1\longmapsto f=x^p+az^p\quad\text{and}\quad q_2\longmapsto g=y^q+bz^q.
$$
The complex obtained by extending this resolution along $\varphi$ is isomorphic to the minimal graded resolution in Theorem~\ref{thm:dichotomy}(2).
\end{proposition}
\begin{proof}
Conjugation is $R$-linear and satisfies $\overline{\bar q_1 u+\bar q_2 v}=\bar u q_1+\bar v q_2$. It therefore identifies $M_2$ with the underlying $R$-module of $\mathcal M$. Under this identification, the resolution above resolves $\mathcal M$ as an $R$-module. For $u,v\in\mathcal A$, write $N(u)=u\bar u=\bar u u\in R$ and set
$$
B(u,v)=\begin{pmatrix}v\bar u&-N(u)\\N(v)&-u\bar v\end{pmatrix}.
$$
For $n=2$, the columns of the syzygy blocks in \cite[Theorem~3.1(1)]{AdamsLoustaunau} are the right multiples of $(q_1\bar q_2,-N(q_1))^t$ and $(N(q_2),-q_2\bar q_1)^t$ by the basis elements $e=1,i,j,k$. Conjugation followed by transposition takes each column to the left multiple by $\bar e$ of the corresponding row of $B(q_1,q_2)$. Thus, with $C=(q_1,q_2)^t$ and $L=\ker\rho_C\subseteq\mathcal A(-1)^2$, we have $\operatorname{im}\rho_{B(q_1,q_2)}=L$.

In the ordinary grading, the surjection $R(-3)^8\to L$ in the resolution of $M_2$ is an isomorphism in degree $3$, since its kernel is isomorphic to $R(-4)^4$. The eight rows obtained by conjugation and transposition span $L_3$ and therefore form a basis. The change of basis in degree $3$ identifies $\rho_{B(q_1,q_2)}\colon\mathcal A(-3)^2\to L$ with that surjection. In particular, $\ker\rho_{B(q_1,q_2)}\cong R(-4)^4$. Put $T=(-\bar q_2,\bar q_1)$. Since the norms are central,
$$
TB(q_1,q_2)=
\bigl(-N(q_2)\bar q_1+\bar q_1N(q_2),\
\bar q_2N(q_1)-N(q_1)\bar q_2\bigr)=(0,0).
$$
The rows $T,iT,jT,kT$ are linearly independent over $\mathbb R$, since $cT=0$ for $c\in\mathbb H$ implies $c=0$. They form a basis of the degree $4$ component of the kernel and generate it over $R$.

Hence $\ker\rho_{B(q_1,q_2)}=\mathcal A T$. Also, $\rho_T$ is injective since $\mathcal A$ is a domain and $T\ne0$. The entries of $C,B(q_1,q_2),T$ have positive degree, so the free resolution over $\mathcal A$ is minimal in the ordinary grading. Its restriction to $R$ is minimal as well, since $\mathcal A_+=R_+\mathcal A$. By uniqueness, this restriction is isomorphic to the resolution of $M_2$ recalled above.

Now assign $t_{1\ell}$ degree $p$ and $t_{2\ell}$ degree $q$. The entries of $C$ have degrees $p,q$, those of $B(q_1,q_2)$ have degrees $p+q,2p,2q,p+q$ in row order, and those of $T$ have degrees $q,p$. The resulting graded free resolution is
$$
0\to\mathcal A(-2p-2q)\xrightarrow{\rho_T}
\mathcal A(-2p-q)\oplus\mathcal A(-p-2q)
\xrightarrow{\rho_{B(q_1,q_2)}}\mathcal A(-p)\oplus\mathcal A(-q)
\xrightarrow{\rho_C}\mathcal A\to\mathcal M\to0.
$$

For $a=a_0+a_1i+a_2j+a_3k$ and $b=b_0+b_1i+b_2j+b_3k$, the assignments $\varphi(t_{1\ell})=\delta_{\ell0}x^p+a_\ell z^p$ and $\varphi(t_{2\ell})=\delta_{\ell0}y^q+b_\ell z^q$ define the graded map $\varphi$, where $\delta_{\ell0}=1$ if $\ell=0$ and $\delta_{\ell0}=0$ otherwise. Since its images on $R$ have real coefficients, $\varphi$ commutes with conjugation. Regard $A$ as a right $\mathcal A$-module by $c\cdot u=c\varphi(u)$. For each shift $d$, the map $A\otimes_{\mathcal A}\mathcal A(-d)\to A(-d)$ defined by $c\otimes u\mapsto c\varphi(u)$ is a graded isomorphism. Under these isomorphisms, tensoring replaces each matrix entry $u$ by $\varphi(u)$.

Since $\varphi(q_1)=f$, $\varphi(q_2)=g$, and $\varphi$ commutes with conjugation, applying $A\otimes_{\mathcal A}-$ to the resolution above yields the complex
$$
\begin{aligned}
0\longrightarrow A(-2p-2q)&\xrightarrow{\rho_{R_0}}A(-2p-q)\oplus A(-p-2q)\\
&\xrightarrow{\rho_{B(f,g)}}A(-p)\oplus A(-q)
\xrightarrow{\rho_{(f,g)^t}}A\longrightarrow M\longrightarrow0,
\end{aligned}
$$
where $R_0=(-\bar g,\bar f)$ and right exactness identifies $A\otimes_{\mathcal A}\mathcal M$ with $M=A/(Af+Ag)$. For $u,v,w\in A$, the three maps are defined by $\rho_{(f,g)^t}(u,v)=uf+vg$, $\rho_{B(f,g)}(u,v)=(u,v)B(f,g)$, and $\rho_{R_0}(w)=(-w\bar g,w\bar f)$.

To compare the matrices, put $h=z^{p+q}$, $S=-\kappa^{-1}I_2$, and $\lambda=-\kappa^{-2}$. Let $u=\operatorname{Im}(a)$ and $v=\operatorname{Im}(b)$, so $\kappa=vu-uv$. Since $u^2,v^2\in\mathbb R$, the commutator satisfies
$$
\kappa u+u\kappa=vu^2-u^2v=0
\quad\text{and}\quad
\kappa v+v\kappa=v^2u-uv^2=0.
$$
Hence $\kappa f\kappa^{-1}=\bar f$ and $\kappa g\kappa^{-1}=\bar g$. Since $f+\bar f$ and $g+\bar g$ are central, $\bar f g-g\bar f=f\bar g-\bar g f=gf-fg=\kappa h$. Using $P_2=hI_2+P_1P_3$, we obtain
$$
\begin{aligned}
S^{-1}P_2
&=-\kappa hI_2+\begin{pmatrix}\bar f\\\bar g\end{pmatrix}(g,-f)
=\begin{pmatrix}g\bar f&-N(f)\\N(g)&-f\bar g\end{pmatrix}=B(f,g),\\
P_3 S
&=(\kappa^{-1}g\kappa^{-1},-\kappa^{-1}f\kappa^{-1})
=\kappa^{-2}(\bar g,-\bar f)=\lambda R_0.
\end{aligned}
$$
By $SB(f,g)=P_2$ and $P_3S=\lambda R_0$, the maps $\mathrm{id},\mathrm{id},\rho_S,\rho_\lambda$ in homological degrees $0,1,2,3$, respectively, define a chain isomorphism from the resolution in Theorem~\ref{thm:dichotomy}(2) to the specialized complex. The maps have degree zero because $S$ is diagonal with entries in $\mathbb H$ and $\lambda\in\mathbb H^\times$. This proves that the specialized complex is a minimal graded free resolution of $M$.
\end{proof}

\section{Betti tables under extension and restriction of scalars}\label{sec:center}
Let $D$ be a division ring and let $K\subseteq Z(D)$ be a subfield with $\delta=\dim_KD<\infty$. For $n\ge1$, put $S=K[x_1,\dots,x_n]$ and $A=D[x_1,\dots,x_n]$, with central indeterminates of degree one. Write $\res{M}$ for the restriction of a left $A$-module $M$ to $S$.

\begin{theorem}\label{thm:restriction}
With this notation, the following hold.
\begin{enumerate}[label=\textup{(\arabic*)}]
\item If $M$ is a finitely generated graded left $A$-module, then $\res{M}$ is a finitely generated graded $S$-module and $\beta^S_{i,j}(\res{M})=\delta\beta^A_{i,j}(M)$ for all $i,j$. If $M\ne0$, then $\pd_S(\res{M})=\pd_A(M)$ and $\reg_S(\res{M})=\reg_A(M)$. Moreover, $\dim_K\res{M}=\delta\dim_DM$, so $M$ is finite-dimensional over $D$ if and only if $\res{M}$ is finite-dimensional over $K$.
\item If $N$ is a finitely generated graded $S$-module, then $A\otimes_SN$ is a finitely generated graded left $A$-module and $\beta^A_{i,j}(A\otimes_SN)=\beta^S_{i,j}(N)$ for all $i,j$. If $N\ne0$, then $A\otimes_SN\ne0$, $\pd_A(A\otimes_SN)=\pd_S(N)$, and $\reg_A(A\otimes_SN)=\reg_S(N)$. Moreover, $\dim_D(A\otimes_SN)=\dim_KN$, so $N$ is finite-dimensional over $K$ if and only if $A\otimes_SN$ is finite-dimensional over $D$.
\end{enumerate}
\end{theorem}

\begin{proof}
The map $D\otimes_KS\to A$ defined by $d\otimes s\mapsto ds$ is an isomorphism of graded $K$-algebras, since $S\subseteq Z(A)$ and the monomials in $x_1,\dots,x_n$ form bases over $K$ and $D$, respectively. If $d_1,\dots,d_\delta$ is a $K$-basis of $D$, then $A=\bigoplus_{r=1}^\delta Sd_r=\bigoplus_{r=1}^\delta d_rS$. In particular, $A$ is free of rank $\delta$ on either side over $S$, and $A_+=S_+A$.

\textup{(1)} Let $F_\bullet\to M$ be a minimal graded free resolution. Restriction to $S$ preserves exactness and the grading, and each $A(-j)$ restricts to $S(-j)^\delta$. In particular, $\res{F_0}\to\res{M}$ shows that $\res{M}$ is finitely generated. By minimality, $\partial_i(F_i)\subseteq A_+F_{i-1}=S_+F_{i-1}$ for $i\ge1$. Thus the restricted resolution is minimal and $\beta^S_{i,j}(\res{M})=\delta\beta^A_{i,j}(M)$. If $(m_\alpha)$ is a left $D$-basis of $M$, then $(d_r m_\alpha)_{r,\alpha}$ is a $K$-basis of $\res{M}$, so $\dim_K\res{M}=\delta\dim_D M$.

\textup{(2)} Let $G_\bullet\to N$ be a minimal graded free resolution over $S$. Since $A$ is free as a right $S$-module, tensoring with $A$ over $S$ preserves exactness. The summand $S(-j)$ extends to $A(-j)$, and the matrices of the differentials have the same entries, now considered in $S_+\subseteq A_+$. It follows that $A\otimes_SG_\bullet$ is a minimal graded free resolution and $\beta^A_{i,j}(A\otimes_SN)=\beta^S_{i,j}(N)$. The surjection $A\otimes_SG_0\to A\otimes_SN$ shows that $A\otimes_SN$ is finitely generated. As graded left $D$-modules, $A\otimes_SN\cong (D\otimes_KS)\otimes_SN\cong D\otimes_KN$. A $K$-basis of $N$ induces a $D$-basis of $D\otimes_K N$, so $\dim_D(A\otimes_SN)=\dim_K N$.

Since $\dim_K\res{M}=\delta\dim_D M$ and $0<\delta<\infty$, the module $\res{M}$ is finite-dimensional over $K$ if and only if $M$ is finite-dimensional over $D$. Likewise, the equality $\dim_D(A\otimes_SN)=\dim_K N$ shows that $A\otimes_SN$ is finite-dimensional over $D$ if and only if $N$ is finite-dimensional over $K$, and that $A\otimes_SN\ne0$ whenever $N\ne0$.

The Betti equalities proved in \textup{(1)} and \textup{(2)}, with $\delta>0$, imply
$$
\begin{aligned}
\beta^S_{i,j}(\res{M})\ne0
&\iff \beta^A_{i,j}(M)\ne0,\\
\beta^A_{i,j}(A\otimes_SN)\ne0
&\iff \beta^S_{i,j}(N)\ne0
\end{aligned}
$$
for every $i,j$. For each nonzero module, $\pd$ is the maximum of $i$ over the pairs with $\beta_{i,j}\ne0$, and $\reg$ is the maximum of $j-i$. The two equivalences therefore prove the equalities of projective dimension and regularity.
\end{proof}

Let $R=D'[x_1,\dots,x_n]$ be a polynomial ring with central indeterminates over a division ring $D'$. Write $\cB(R)$ for the set of Betti tables of nonzero finitely generated graded left $R$-modules of finite length. For two such modules $M,N$, the module $M\oplus N$ has finite length, and the direct sum of their minimal graded free resolutions is minimal. Hence $\beta^R(M\oplus N)=\beta^R(M)+\beta^R(N)$, so $\cB(R)$ is a semigroup under addition. Following the terminology over a field in \cite[Definitions~1.1--1.2, p.~343]{ErmanSemigroup}, define the rational cone generated by $\cB(R)$ by
$$
\cC(R)=\left\{\sum_{\nu=1}^{r}\lambda_\nu\gamma_\nu\ \middle|\
r\in\mathbb Z_{\ge1},\ \lambda_\nu\in\mathbb Q_{\ge0},\ \gamma_\nu\in\cB(R)
\right\}.
$$

\begin{corollary}\label{cor:betti-cones}
With $S$ and $A$ as in Theorem~\ref{thm:restriction}, we have $\cB(S)\subseteq\cB(A)$ and $\delta\cB(A)\subseteq\cB(S)$. In particular, $\cC(S)=\cC(A)$. If $n=3$ and $D$ is noncommutative, then $\cB(S)\subsetneq\cB(A)$.
\end{corollary}

\begin{proof}
Let $N$ be a nonzero finitely generated graded $S$-module of finite length. By Proposition~\ref{prop:finite-length} and Theorem~\ref{thm:restriction}\textup{(2)},
$$
0<\dim_D(A\otimes_SN)=\dim_K N<\infty
\quad\text{and}\quad
\beta^A(A\otimes_SN)=\beta^S(N).
$$
Thus $A\otimes_SN$ has finite length, and $\beta^S(N)\in\cB(A)$. Hence $\cB(S)\subseteq\cB(A)$.

For a nonzero finitely generated graded $A$-module $M$ of finite length, Theorem~\ref{thm:restriction}\textup{(1)} implies
$$
0<\dim_K\res{M}=\delta\dim_D M<\infty
\quad\text{and}\quad
\beta^S(\res{M})=\delta\beta^A(M).
$$
Again by Proposition~\ref{prop:finite-length}, $\res{M}$ has finite length. Hence $\delta\beta^A(M)\in\cB(S)$, so $\delta\cB(A)\subseteq\cB(S)$.

The inclusion $\cB(S)\subseteq\cB(A)$ implies $\cC(S)\subseteq\cC(A)$. For $\gamma\in\cB(A)$, we have $\delta\gamma\in\cB(S)$ and $\gamma=\delta^{-1}(\delta\gamma)\in\cC(S)$, since $\delta^{-1}\in\mathbb Q_{>0}$. Thus $\cC(A)\subseteq\cC(S)$, and the cones are equal.

Suppose $n=3$ and $D$ is noncommutative. Choose $a,b\in D$ with $ba-ab\ne0$ and put $M=A/(A(x+az)+A(y+bz))$. By Theorem~\ref{thm:dichotomy}\textup{(2)}, $\dim_D M=2$ and $\beta^A(M)=D_1$. Proposition~\ref{prop:finite-length} and Corollary~\ref{prop:notcomm} therefore imply $D_1\in\cB(A)\setminus\cB(S)$. Together with $\cB(S)\subseteq\cB(A)$, this proves $\cB(S)\subsetneq\cB(A)$.
\end{proof}

\appendix
\section{Computational checks}
\label{app:code}

The supplementary scripts test the matrix identities in Theorem~\ref{thm:dichotomy} on examples with $\kappa=0$ and $\kappa\ne0$ and compare Hilbert functions with counts of standard monomials. They also check examples of the pure case, the calculation for two generators over a field, the quaternion specialization in Proposition~\ref{rem:quaternion-specialization}, and the Betti formulas in Theorem~\ref{thm:restriction}.

We check polynomial and matrix identities over $\mathbb Q$. In internal degree $j$, write $d_{i,j}\colon V_{i,j}\to V_{i-1,j}$ for the differential of the complex of $\mathbb Q$-vector spaces. In the rank tests, the denominators of the matrix entries are prime to $\ell=2^{31}-1$. Write $r_i$ for the rank over $\mathbb Q$ and $\overline r_i$ for the rank after reduction modulo $\ell$. A nonzero minor modulo $\ell$ is nonzero over $\mathbb Q$, so $\overline r_i\le r_i$. Once $d_{i,j}d_{i+1,j}=0$ has been checked over $\mathbb Q$,
$$
0\le\dim_{\mathbb Q}H_i(V_{\bullet,j})
=\dim_{\mathbb Q}V_{i,j}-r_i-r_{i+1}
\le\dim_{\mathbb Q}V_{i,j}-\overline r_i-\overline r_{i+1}.
$$
Thus $\overline r_i+\overline r_{i+1}=\dim_{\mathbb Q}V_{i,j}$ certifies exactness at $V_{i,j}$. At the leftmost term the test is $\overline r_i=\dim_{\mathbb Q}V_{i,j}$, which certifies injectivity.

For the Betti comparisons, let $\overline M$ be the cokernel of a presentation reduced modulo a prime $s$. Let $K_\bullet(x;\overline M)$ be the Koszul complex on the central variables with coefficients in $\overline M$. If $C$ is the coefficient algebra of the presentation over $\mathbb Q$, we compare
$$
\frac{\dim_{\mathbb F_s}H_i(K_\bullet(x;\overline M))_j}{\dim_{\mathbb Q}C}
$$
with the Betti numbers in the text. These comparisons use $s=\ell$ and, for some commutative examples, $s=2,3$. The comparison in each case is over $\mathbb F_s$.

\section*{Declaration of competing interest}
The authors declare that they have no competing interests.

\section*{Data and code availability}
The Python scripts for Appendix~\ref{app:code} are available on Zenodo~\cite{TuKhanhPython}. Section~1 of the accompanying documentation describes how to run them.

\label{pg:bib}

\end{document}